\documentclass[11pt]{article}
\usepackage{graphicx}
\usepackage{color}
\usepackage{amsmath}
\usepackage{amssymb}
\usepackage{amscd}
\usepackage{framed}
\usepackage{bbm}

\usepackage{fancyhdr,a4wide}
\usepackage{amsthm}

\newcommand{\R}{\mathbb{R}}
\newcommand{\inr}[1]{\left\langle #1 \right\rangle}

\newcommand{\E}{\mathbb{E}}

\newcommand{\eps}{\varepsilon}

\newtheorem{Theorem}{Theorem}[section]
\newtheorem{Lemma}[Theorem]{Lemma}

\newtheorem{Definition}[Theorem]{Definition}

\newtheorem{Proposition}[Theorem]{Proposition}

\newtheorem{Corollary}[Theorem]{Corollary}
\newtheorem{Remark}[Theorem]{Remark}
\newtheorem{Example}[Theorem]{Example}
\newtheorem{Assumption}{Assumption}[section]
\newtheorem{Question}[Theorem]{Question}

\numberwithin{equation}{section}

\def \proof {\noindent {\bf Proof.}\ \ }

\def \endproof
{{\mbox{}\nolinebreak\hfill\rule{2mm}{2mm}\par\medbreak}}
\def\IND{\mathbbm{1}}

\newcommand{\PROB}{\mathbb{P}}

\def\IND{\mathbbm{1}}

\begin{document}
\title{Stable recovery and the coordinate small-ball behaviour of random vectors}

\author{Shahar Mendelson \thanks{LPSM, Sorbonne University, Paris, France, and MSI, The Australian National University, Canberra, Australia; \ shahar.mendelson@upmc.fr.}
\and
Grigoris Paouris \thanks{Department of Mathematics, Texas A$\&$M University, College Station, TX 77843-3368, U.S.A.; \ grigoris@math.tamu.edu. Supported by NSF  grant DMS-1812240.}}

\date{}

\maketitle

\begin{abstract}
Recovery procedures in various application in Data Science are based on \emph{stable point separation}. In its simplest form, stable point separation implies that if $f$ is ``far away" from $0$, and one is given a random sample $(f(Z_i))_{i=1}^m$ where a proportional number of the sample points may be corrupted by noise, that information is still enough to exhibit that $f$ is far from $0$.

Stable point separation is well understood in the context of iid sampling, and to explore it for general sampling methods we introduce a new notion---the \emph{coordinate small-ball} of a random vector $X$. Roughly put, this feature captures the number of ``relatively large coordinates" of $(|\inr{TX,u_i}|)_{i=1}^m$, where $T:\R^n \to \R^m$ is an arbitrary linear operator and $(u_i)_{i=1}^m$ is any fixed orthonormal basis of $\R^m$.

We show that under the bare-minimum assumptions on $X$, and with high probability, many of the values $|\inr{TX,u_i}|$ are at least of the order $\|T\|_{S_2}/\sqrt{m}$. As a result, the ``coordinate structure" of $TX$ exhibits the typical Euclidean norm of $TX$ and does so in a stable way.

One outcome of our analysis is that random sub-sampled convolutions satisfy stable point separation under minimal assumptions on the generating random vector---a fact that was known previously only in a highly restrictive setup, namely, for random vectors with iid subgaussian coordinates.
\end{abstract}

\section{Introduction}
One of the key questions in Data Science is to identify (or at least approximate) an unknown function using partial information. In standard recovery problems the data one receives consists of a finite sample of the unknown function and the sample points are assumed to be independent. The sample is then used to construct a suitable `guess' of the function and the hope is that the guess is a good approximation in some appropriate sense.

Off-hand, the significance of having sample points that are selected independently is not clear. A closer inspection shows that independence has a strong geometric impact: it leads to \emph{point separation}.

\subsection*{Point separation and stable point separation}
To explain what we mean by point separation, let us first consider it in its simplest form, separation of a function from $0$.

Given a function $f$ on a probability space $(\Omega,\mu)$, let $Z$ be distributed according to $\mu$ and consider a sample $Z_1,...,Z_m$, consisting of independent points distributed as $Z$. The sample $(Z_1,...,Z_m)$ naturally endows a random vector $X \in \R^m$, whose coordinates are the given measurements $f(Z_i), \ 1 \leq i \leq m$; that is,
$$
X=\left( f(Z_1),...,f(Z_m)\right).
$$
Any hope of identifying $f$ from the given data vector $X$ is based on the belief that $X$ captures enough features of $f$; for example, that $f$ can be distinguished from $0$ with only $X$ as data. Thus, one has to address the following question:
\begin{Question}
If $f$ is reasonably far away from $0$, when is that fact exhibited by a typical realization of $X$?
\end{Question}

An obvious way of exhibiting separation between $f$ and $0$ is through the Euclidean norm of $X$; specifically, by showing that, with high probability,
\begin{equation} \label{eq:l-2-intro}
\frac{\|X\|_2^2}{m} = \frac{1}{m}\sum_{i=1}^m f^2(Z_i) \geq \kappa \|f\|_{L_2}^2
\end{equation}
for a suitable constant $\kappa$ and for any $m \geq m_0$. Independence proves to be extremely useful in establishing \eqref{eq:l-2-intro}. Indeed, under a weak \emph{small-ball assumption}, that
\begin{equation} \label{eq:small-ball-single}
\mathbb P (|f(Z)| \geq \kappa \|f\|_{L_2}) \geq \rho,
\end{equation}
it is straightforward to verify that with probability at least $1-2\exp(-c \rho m)$,
\begin{equation} \label{eq:single-rv-independence}
|\{i :|f(Z_i)| \geq \kappa \|f\|_{L_2} \} | \geq  \frac{\rho}{2} m.
\end{equation}
Thus, with very high probability, a proportion of the coordinates of $X$ are large---of the order of $\|f\|_{L_2}$.

\vskip0.4cm

While \eqref{eq:single-rv-independence} clearly implies \eqref{eq:l-2-intro} and point separation, it says much more: under the small-ball assumption \eqref{eq:small-ball-single}, \emph{independent sampling leads to stable point separation}: not only is $\|X\|_2$ large, the reason that it is large is because many of its coordinates $|\inr{X,e_i}|$ are nontrivial---making point separation robust to noise. In particular, even if a (small) fraction of the measurements $f(Z_i)$ are corrupted maliciously, the fact that $f$ is far away from $0$ is still exhibited by the corrupted vector.

\vskip0.4cm

In a more geometric language, stable point separation is manifested by the fact that $(\inr{X,e_i})_{i=1}^N$ is a well-spread vector, and obviously this significant additional information does not come for free: stable point separation is much harder to prove than point separation. At the same time, the importance of the notion is clear: intuitively, a sampling method can be useful in statistical recovery problems, where being robust to noise is of the utmost importance, only if it satisfies a \emph{uniform version of stable point separation}. Indeed, at the heart of numerous statistical procedures is the fact that if $F$ is a class of functions, then with high probability, for every $f,h \in F$ that are `far enough'
\begin{equation} \label{eq:intro-uni-stable-sep}
|\{i : |(f-h)(X_i)| \geq \kappa \|f-h\|_{L_2} \} | \geq c(\rho)m,
\end{equation}
which is a uniform version of stable point separation. It allows one to distinguish between any two functions in the given class that are sufficiently far apart using a typical sample, even when a proportional number of the given measurements are corrupted by noise.

Uniform stable point separation has played a central role in the recent progress on some key questions in learning theory and statistics. For example, it has led to the introduction of an optimal learning procedure in \cite{LM1,M1}; to optimal vector mean estimation in \cite{LM2,LM3} and to optimal covariance estimation in \cite{M2,MZ} --- all of which in heavy-tailed situations.

Unfortunately, stable point separation and its uniform counterpart are well understood only for iid sampling, and the downside of iid sampling is that it leads to various computational difficulties. For example, consider a relatively simple recovery problem, where the goal is to identify an unknown $t_0 \in T \subset \R^n$ using linear measurements $(\inr{Z_i,t_0})_{i=1}^m$ and $Z_1,...,Z_m$ are independent copies of the standard Gaussian random vector in $\R^n$. Procedures that aim at recovering $t_0$ are based on vector multiplications with the matrix $\Gamma=\sum_{i=1}^m \inr{Z_i,\cdot}e_i$, but because $\Gamma$ has independent Gaussian rows, vector multiplication is computationally expensive.

To address these and other computational difficulties of a similar nature, other sampling methods are often used  in recovery problems. However, once the iid framework is abandoned, establishing the required point separation/stable point separation becomes a formidable task; in fact, it is often far from obvious that either one of the properties is true when the sample points are not independent.

Motivated by general sampling methods, the question we focus on is as follows: are point separation and stable point separation really the outcome of independence?  Rather informally, the question we study is:
\begin{Question} \label{qu:main-intro}
Given a centred random vector $X \in \R^n$,
\begin{description}
\item{$(a)$} What conditions on $X$ are needed to ensure that for an arbitrary linear operator $T:\R^n \to \R^m$ and with high probability $\|T X\|_2$ is reasonably large?
\item{$(b)$} When is the fact that $\|TX\|_2$ is large exhibited by an arbitrary coordinate structure? In other words, given an arbitrary orthonormal basis $(u_i)_{i=1}^m$, are many of the values $|\inr{TX,u_i}|$ reasonably large?
\end{description}
\end{Question}

Question \ref{qu:main-intro} clearly extends the notions of point separation and stable point separation from the iid setup, where $m=n$ and $X=(f(Z_i))_{i=1}^m$: in the general framework of Question \ref{qu:main-intro} the coordinates of $X$ need not be independent or identically distributed, and $X$ is further distorted by a linear operator $T$.

\vskip0.4cm

Let us illustrate how addressing the two parts of Question \ref{qu:main-intro} can become unpleasant very quickly once independence is left behind. The example we present here is the very popular \emph{random sub-sampled convolutions} scheme, which is used in numerous applications, such as SAR radar imaging, optical imaging, channel estimation, etc. (see \cite{Rom,FocRau} for more details on these and other applications).

\begin{Example} \label{ex:circulant}
Let $\xi$ be an isotropic random vector in $\R^n$ (that is, $\xi$ is centred and for every $t \in \R^n$, $\E\inr{\xi,t}^2 = \|t\|_2^2$). Fix $a \in \R^n$ and let $W=a \circledast \xi$ be the discrete convolution of $a$ and $\xi$; i.e., if $j \ominus i = j - i \ {\rm mod} \ n$ and $\tau_i$ is the shift operator defined by $(\tau_i x)=(x_{j \ominus i})_{j=1}^n$, then
$$
a \circledast \xi = (\inr{a, \tau_i \xi})_{i=1}^n.
$$

The measurements of the vector $a$ one receives come from a selection of a random subset of the coordinates of $a \circledast \xi$: let $\delta_1,...,\delta_n$ be independent, $\{0,1\}$-valued random variables with mean $\delta$; set $I=\{i : \delta_i=1\}$; and define $Z=(a \circledast \xi)_{i \in I}$.

Note that typically  $|I| \sim \delta n$ and $\E \|Z\|_2^2=\delta n \|a\|_2^2$. Therefore, this sampling method exhibits point separation of $a$ and $0$ if, with high probability,
\begin{equation} \label{eq:part-conv-intro}
\frac{1}{\delta n} \sum_{i \in I} Z_i^2 \geq c \|a\|_2^2
\end{equation}
for a suitable constant $c$ (that should be independent of $a$ and $\delta$). And, it exhibits stable point separation of $a$ and $0$ with respect to the standard basis $(e_i)_{i=1}^n$ if with high probability,
\begin{equation} \label{eq:part-conv-intro-stable}
\left|\left\{i \in I : |\inr{Z,e_i}| \geq c^\prime \|a\|_2 \right\} \right| \geq c^{\prime \prime} \delta n
\end{equation}
for suitable constants $c^\prime$ and $c^{\prime \prime}$. In particular, \eqref{eq:part-conv-intro-stable} means that the large Euclidean norm of $Z$ is exhibited by the fact that many of the coordinates of $Z$ (with respect to the standard basis) are large.

The advantage in using the random sub-sampled convolution scheme is that recovery can be carried out rather efficiently: for instance, unlike matrices with iid rows, there is a fast matrix--vector multiplication algorithm for partial circulant matrices (see, e.g., \cite{GolLoa}).

\end{Example}

Clearly, identifying when, or even if, \eqref{eq:part-conv-intro} and \eqref{eq:part-conv-intro-stable} are true is considerably harder than establishing \eqref{eq:l-2-intro} and \eqref{eq:single-rv-independence}. And if they are, it has nothing to do with independence.

\vskip0.4cm

A wildly optimistic conjecture is that both parts of Question \ref{qu:main-intro} are (almost) universally true under minimal assumptions on $X$. And deferring an accurate definition of what is meant by ``reasonably large", the main result of this article is that this wildly optimistic conjecture is, in fact, true:
\begin{framed}
\begin{description}
\item{$\bullet$}  Under the bare-minimum assumptions on $X$, for an arbitrary linear operator $T$, $TX$ has a large Euclidean norm; moreover, that norm is exhibited by many large coordinates \emph{with respect to an arbitrary orthonormal basis}.

\item{$\bullet$} Both facts hold with high probability and are simply generic properties of $X$ that have nothing to do with independence, nor with concentration of measure.

\item{$\bullet$} In particular, almost any random vector $TX$ exhibits both point separation and stable point separation with respect to an arbitrary orthonormal basis.
\end{description}
\end{framed}

We show in what follows that the reason why both parts of Question \ref{qu:main-intro} are universally true is a \emph{small-ball assumption} which we now describe.

\subsection{The Small-ball assumption}
To have some intuition on the sort of quantitative answers to Question \ref{qu:main-intro} one can hope for, assume for the time being that $X$ is isotropic. Let $F \subset \R^n$ be a subspace of dimension $k$, and set $P_F$ to be the orthogonal projection onto $F$. Thus, $\E\|P_F X\|_2^2 = k$ and at least intuitively, saying that  $P_F X$ has a ``reasonably large Euclidean norm" can be taken to mean that $\|P_F X\|_2 \geq \eps \sqrt{k}$ for some $0<\eps<1$. Moreover, a ``reasonably large coordinate" of such a $k$-dimensional vector should be at least of the order of $(\E \|P_F X\|_2^2)^{\frac{1}{2}}/\sqrt{k}=1$.

Following the same path with a general linear operator $T:\R^n \to \R^m$ instead of $P_F$, the intuitive notion of being relatively large is that $\|TX\|_2$ is at least $\eps \|T\|_{S_2}=\eps(\E \|TX\|_2^2)^{\frac{1}{2}}$, where $\|T\|_{S_2}$ denotes the Hilbert-Schmidt norm of $T$; and given an orthonormal basis $(u_i)_{i=1}^m$, a ``large coordinate'' of $TX$ satisfies that $|\inr{TX,u_i}| \gtrsim \|T\|_{S_2}/\sqrt{m}$.

\vskip0.4cm

Once the two notions are agreed upon, the answers to the two parts of Question \ref{qu:main-intro} are given in the form of \emph{small-ball estimates}, that is, upper bounds on
$$
\PROB(\|TX\|_2 \leq y) \ \ \ \ \ {\rm for} \ \ \ \ \ 0<y \leq \|T\|_{S_2}
$$
and \emph{coordinate small-ball estimates} which are upper bounds on
$$
\PROB \left( \left | \{i : |\inr{TX,u_i}| \leq y \} \right| \geq \ell \right) \ \ \ \ \ {\rm for} \ \ \ \ \ 0 < y \leq \|T\|_{S_2}/\sqrt{m}.
$$

Let us emphasize a fact, which at first glance, may be surprising:

\vskip0.4cm
\begin{framed}
Small-ball estimates and coordinate small-ball estimates have nothing to do with concentration.
\end{framed}

Indeed, although the notions of small-ball estimates and coordinate small-ball estimates may seem to be related to \emph{concentration of measure}, they are actually based on a totally different phenomenon that has nothing to do with the way the random variable $\|TX\|$ concentrates around its mean $\E\|TX\|$---no matter what norm $\| \ \|$ is considered.

There are several reasons for that: firstly, two-sided concentration estimates of the form
$$
\PROB \left( \left| \|TX\| - \E\|TX\| \right| \geq y \right)
$$
are a combination of the upper estimate---that with high probability, $\|TX\| \leq \E\|TX\|+y$,  and the lower one, that $\|TX\| \geq \E\|TX\|-y$. By now it is well understood (see, for example, the discussion in \cite{Men}) that the two estimates are totally different and are caused by unrelated features of the random vector $X$. Moreover, the upper tail is almost always the bottleneck in the two-sided estimate, while our interests lie in the lower one. Secondly, the scale one is interested in when studying the small-ball behaviour of $TX$ corresponds to the lower tail with the choice of $y=(1-s)\E\|TX\|$, for $s$ close to $0$. That is very different from the lower tail at the `concentration scale' of $y=s\E\|TX\|$ for $s$ close to $0$.

\vskip0.4cm

\begin{Remark}
As we explain in what follows, the behaviour of $\PROB(\|TX\|_2 \leq y)$ is more subtle than what this intuitive description may lead one to believe. In fact, $\PROB(\|TX\|_2 \leq \eps \|T\|_{S_2})$  exhibits multiple phase transitions at different levels of $\eps$ in the small-ball regime.
\end{Remark}

The minimal assumption that is required for establishing small-ball and coordinate small-ball estimates is as follows:

\begin{Assumption} \label{ass:small-ball}
The random vector $X$ satisfies a \emph{small ball assumption} (denoted from here on by {\it SBA}) with constant ${\cal{L}}$ if for every $ 1\leq k \leq n-1$,  every $ k$ dimensional subspace $F$, every $z\in \R^{n}$ and every $\eps>0$,
\begin{equation} \label{sb-asum-1-1}
\mathbb P \left ( \| P_{F} X- z\|_{2} \leq \eps \sqrt{k} \right) \leq ( {\cal{L}} \eps)^{k},
\end{equation}
where $P_F$ is the orthogonal projection onto the subspace $F$.
\end{Assumption}
It is straightforward to verify (see, e.g.,  Proposition 2.2 in \cite{RV}) that $X$ satisfies the {\it SBA} with constant ${\cal L}$ if and only if for every $1 \leq k \leq n$, the densities of all $k$-dimensional marginals of $X$ are bounded by $ {\cal{L}}^{k}$ (assuming, of course, that $X$ has a density, and in which case $f_{X}$ denotes that density).

There are numerous examples of generic random vectors that satisfy the {\it SBA} with an absolute constant; among them are vectors with iid coordinates that have a bounded density  (\cite{RV} and \cite {LPP} for the optimal constant)  as well as various log-concave random vectors\footnote{Recall that $X$ is log-concave if it has a density $ f_{X}$ that satisfies that for every $ x, y$ in the support of $ f_{X}$ and every $0 \leq \lambda \leq 1$,
$f_{X} ((1- \lambda) x+ \lambda y )\geq f_{X}^{(1-\lambda)} (x) f_{X}^{\lambda} (y)$.}. For more details see Appendix \ref{app:log-concave0}, where we list several examples of generic log-concave random vectors that satisfy Assumption \ref{ass:small-ball}.

\vskip0.4cm

Although Assumption \ref{ass:small-ball} requires that $X$ has a density, this is not essential and our main results remain true even under the following weaker assumption.

 \begin{Assumption} \label{ass:w-small-ball}
Let ${\cal L}$ and $\theta$ be such that ${\cal L}\theta<1$. The random vector $X$ satisfies the \emph{weak small-ball assumption} (denoted from here on by {\it wSBA}) with constants $\theta$ and  ${\cal{L}}$ if  for every $ 1\leq k \leq n-1$,  every $ k$ dimensional subspace $F$, and every $z\in \R^{n}$,
\begin{equation} \label{sb-asum-1-w}
\mathbb P \left ( \| P_{F} X- z\|_{2} \leq \theta \sqrt{k} \right) \leq ( {\cal{L}} \theta)^{k}.
\end{equation}
\end{Assumption}

Clearly, if $X$ satisfies the {\it SBA} with constant ${\cal{L}}$ then it satisfies the {\it wSBA} with constants $ \theta$ and ${\cal{L}}$ for every $ \theta>0$. Moreover, it follows from \cite{RV} that if $ X$ has independent coordinates, and if each coordinate satisfies the {\it wSBA} with constants $ \theta$ and ${\cal{L}}$, then $ X$ satisfies the {\it wSBA} with constants $C \theta$ and ${\cal{L}}$, where $ C>0$ is an absolute constant.

\begin{Remark}
Most of the results presented in what follows hold under the {\it wSBA}. However, to simplify the presentation only one result is proved under that assumption---the coordinate small-ball estimate (Theorem \ref{th2-main-s2}); the other results are formulated using the {\it SBA} which leads to a proof that is less involved.
\end{Remark}

\vskip0.4cm

Before we formulate the main results, let us mention one of their outcomes: a stable point separation bound for the random sub-sampled convolutions scheme.

\subsubsection*{Example \ref{ex:circulant} revisited}

As it happens, the existing state of the art on point separation/stable point separation of the random sub-sampled convolutions scheme can be improved dramatically, as existing estimates are based on severe restrictions on the random vector $\xi$. The reason for those restrictions is a wasteful method of proof, as is explained in Section \ref{sec:wrong-way}, and which leads to the following:

\begin{Theorem} \label{thm:SB-pt-circulant-old} \cite{MRW}
For every constant $L \geq 1$ there exist constants $c_0,c_1,c_2,c_3$ and $c_4$ that depend only on $L$ for which the following holds. Let $x$ be a mean-zero, variance one, $L$-subgaussian random variable\footnote{A centred random variable $x$ is $L$-subgaussian if for every $p \geq 2$, $\|x\|_{L_p} \leq L \sqrt{p}\|x\|_{L_2}$.}, and set $\xi=(x_i)_{i=1}^n$, i.e., a vector whose coordinates are independent copies of $x$. Let $s \leq c_0 n/\log^4 n$ and consider $a \in S^{n-1}$ that is $s$-sparse with respect to the standard basis $(e_i)_{i=1}^n$. Then with probability at least $1-2\exp(-c_1\min\{n/s,\delta n\})$ with respect to both $\xi$ and $(\delta_i)_{i=1}^n$,
\begin{equation} \label{eq:CSB-circulant-old}
\sum_{i \in I} \inr{a \circledast \xi,e_i}^2 \geq c_2\delta n \ \ {\rm and} \ \    \left|\{i \in I : |\inr{a \circledast \xi,e_i}| \geq c_3\}\right| \geq c_4 \delta n.
\end{equation}
\end{Theorem}

We show that one can replace the wasteful parts of Theorem \ref{thm:SB-pt-circulant-old}, leading to a sharp point separation and stable point separation that hold as long as $\xi$ satisfies the {\it SBA}, and with a much better probability estimate.

To formulate this fact, let ${\cal F}$ denote the un-normalized discrete Fourier matrix in $\R^n$ (in particular, ${\cal F}$ is a matrix whose entries are either $-1$ or $1$), and for $a \in \R^n$, let $\hat{a}={\cal F}a/\sqrt{n}$ be the normalized discrete Fourier transform of $a$.

\begin{Theorem} \label{thm:circulant-small-ball}
Let $\xi$ satisfies the {\it SBA} with constant ${\cal L}$ and consider $a \in S^{d-1}$ that is $s$-sparse for $s \leq c_0n/\log n$. Then for any $0<\eps<1$ and $q>2$, with probability at least
$$
1-(c_1 {\cal L}\eps)^{\frac{c_2}{\|\hat{a} \|_{q}^{2q/(q-2)}}}-\exp(-c_3 \delta n),
$$
$$
\sum_{i \in I} \inr{a \circledast \xi,e_i}^2 \geq c_4 \eps^2 \delta n \ \ {\rm and} \ \  \left| \{i \in I : |\inr{a \circledast \xi,e_i}| \geq \eps \} \right| \geq c_5 \delta n;
$$
here, $c_0,c_1$ and $c_2$ are constants that depend on $q$ and $c_3,c_4, c_5$ are absolute constants.
\end{Theorem}

The differences between Theorem \ref{thm:SB-pt-circulant-old} and Theorem \ref{thm:circulant-small-ball} are substantial. Firstly, the estimate in Theorem \ref{thm:circulant-small-ball} holds for a random vector that satisfies the SBA rather than only for vectors that have iid subgaussian coordinates. Secondly, note that for any $a \in S^{n-1}$ and any $q>2$, $1/\|\hat{a}\|_q^{2q/(q-2)} \geq 1/\|\hat{a}\|_\infty^2$; and for any $a \in S^{n-1}$ that is $s$-sparse, $1/\|\hat{a}\|_\infty^2 \geq n/s$. Thus, the probability estimate in Theorem \ref{thm:circulant-small-ball} is always better than in Theorem \ref{thm:SB-pt-circulant-old}, and often the gap between the two is significant.

\begin{Remark}
It is possible to prove a version of Theorem \ref{thm:circulant-small-ball} for $X$ that satisfies the {\it wSBA}, but for the sake of a simpler presentation we shall not do that.
\end{Remark}

\subsection{Small-ball estimates}
If one wants to highlight the crucial (and rather remarkable) feature of the small-ball estimate presented here, it is the following:
\begin{framed}
A Gaussian random vector is not the best case; actually, it is the worst one.
\end{framed}

To explain what we mean by this, let $X$ be a random vector taking values in $\R^n$ and for now assume that it satisfies the {\it SBA} with constant ${\cal L}$. Set $m \leq n$ and let $T:\R^n \to \R^m$ be a linear operator of full rank. Without loss of generality one may assume that $T$ actually maps $\R^n$ into $\R^n$ and denote by $ s_{1}, \cdots , s_{m}$ its nonzero singular values.

Recall that the $p$-Schatten norm of $T$ is
$$
\|T\|_{S_p} = \left(\sum_{i=1}^m s_i^p \right)^{1/p},
$$
and following \cite{NY}, for $2<q \leq \infty$, let
$$
{\rm srank}_q (T) = \left(\frac{\|T\|_{S_2}}{\|T\|_{S_q}}\right)^{\frac{2q}{q-2}}
$$
be the $q$-{\it stable rank} of $T$. Clearly ${\rm srank}_q (T)  \leq m= {\rm rank}(T)$ and the case $q=\infty$ corresponds to the standard notion of the stable rank, i.e.,
$$
{\rm srank}(T) = \left(\frac{\|T\|_{S_2}}{\|T\|_{S_\infty}}\right)^2.
$$

The current state of the art as far as small-ball estimates are concerned is due to Rudelson and Vershynin:

\begin{Theorem} \label{thm:RV} \cite{RV}
There are absolute constants $c_0$ and $c_1$ for which the following holds. If $X$ satisfies the {\it SBA} with constant $ {\cal{L}}$ then for any $\eps>0$,
\begin{equation} \label{SB-1-2}
\mathbb P \left( \| T X\|_{2} \leq \varepsilon \| T \|_{S_{2}}\right) \leq ( c {\cal{L}} \varepsilon)^{ c_{1} {\rm srank}(T) }.
\end{equation}
\end{Theorem}

\begin{Remark}
Although Theorem \ref{thm:RV} is not stated explicitly in \cite{RV}, it follows from the analysis presented there in a straightforward way. Previous estimates of the same flavour have been derived for a centred random vector $X$ that has independent subgaussian entries in \cite{LMOT} and for an $X$ that is isotropic, log-concave and subgaussian in \cite{P-sb}.
\end{Remark}

One instance in which Theorem \ref{thm:RV} can be applied is when $T$ is an orthogonal projection of rank $k$ (and in which case, $ {\rm srank}(T)= k$). On the other hand, it is straightforward to verify that if \eqref{SB-1-2} holds for any such orthogonal projection then $X$ satisfies the {\it SBA} (though perhaps with a slightly different constant). Despite this equivalence, Theorem \ref{thm:RV} is far from optimal --- because of the loose probability estimate; it does not ``see" phase transitions that occur as $\eps$ decreases.

\vskip0.4cm

In contrast to Theorem \ref{thm:RV}, our first main result is a comparison theorem which shows that the worst random vector in the context of small-ball estimates is actually the standard Gaussian. Then, in Corollary \ref{cor:SB-intro} and Theorem \ref{thm:small-ball-gaussian} one uses the Gaussian case to establish the right probability estimate at every scale.
\begin{Theorem} \label{main-compl-sb}
 Let $ X$ be an $n$-dimensional random vector that satisfies the {\it SBA} with constant ${\cal{L}}$, let $T: \mathbb R^{n}\rightarrow \mathbb R^{m}$ be a linear map and set $G$ to be the standard Gaussian vector in $\R^n$. Then for every $ 1\leq k < m = {\rm rank}(T)$,
 \begin{equation} \label{main-compl-sb-1}
 \E\left(\|T X\|_{2}^{-k}\right)   \leq \E \left(\|T G/(\sqrt{2\pi}{\cal{L}})\|_{2}^{-k}\right).
 \end{equation}
 \end{Theorem}

The connection between small-ball estimates and negative moments is an immediate corollary of Markov's inequality, which, combined with Theorem \ref{main-compl-sb}, implies that
$$
\mathbb P \left( \|TX\|_2 \leq \eps \right) = \mathbb P \left( \|TX\|_2^{-k} \geq \eps^{-k} \right) \leq \eps^k \E\left(\|TX\|_2^{-k}\right) \leq  \eps^k \E\left(\|T G/(\sqrt{2\pi}{\cal{L}})\|_{2}^{-k}\right).
$$

\begin{Remark} \label{rem:sharp}
It is natural to ask whether Theorem \ref{main-compl-sb} is sharp, as potentially there could be a significant gap between $(\E\|T X\|_{2}^{-k})^{-\frac{1}{k}}$ and $( \E \|T G\|_{2}^{-k})^{-\frac{1}{k}}$. However, the two happen to be equivalent for any centred log-concave measure (up to the {\it SBA} constant $ {\cal{L}}$). Indeed, one can show that there is an absolute constant $c_0$ such that for any centred log-concave random vector $X$ and $1 \leq k \leq {\rm rank}(T)$,
\begin{equation} \label{main-compl-sb-1-log}
 \left( \E \|T X\|_{2}^{-k} \right)^{-\frac{1}{k}}  \leq c_{0} \left( \E \|T G\|_{2}^{-k} \right)^{-\frac{1}{k}}.
\end{equation}
A sketch of the proof of this fact is presented in Appendix \ref{app:log-concave}.
\end{Remark}

Theorem \ref{main-compl-sb} is a clear indication that the small-ball behaviour of a random vector has nothing to do with concentration or with tail estimates: to a certain extent, concentration exhibited by Gaussian vectors is the best one can hope for, but when it comes to small-ball estimates the situation is the complete opposite. Moreover, thanks to the lower bound from Theorem \ref{main-compl-sb}, the worst case scenario is actually very good and can be controlled. Indeed, to complement Theorem \ref{main-compl-sb} one may estimate the negative moments of $\|TG\|_2$---which requires the following definition:

\begin{Definition} \label{def:a-k} Let $ T:\mathbb R^{n} \rightarrow \mathbb R^{m} $ and $ {\rm rank} ( T) = m $.
For $1 \leq k \leq m-1$ set,
\begin{equation} \label{def-an}
 a_{k} ( T ) = \left( \int_{{\cal G}_{m,k}} {\rm det}^{-\frac{1}{2}}[ (P_{F} T) ( P_{F}T )^{\ast}] dF \right)^{-\frac{1}{ k}}
\end{equation}
where $P_F$ is the orthogonal projection onto the subspace $F$ and the integration takes place on the Grassmannian ${\cal G}_{m,k}$ with respect to the Haar measure. Also for $k=m$ put
$$
a_{m} (T)= {\rm det}^{\frac{1}{2m}}(T  T ^{\ast} ).
$$
\end{Definition}

It is straightforward to verify that $a_k(T)$ has strong ties to the negative moments of $\|TG\|_2$. Indeed, as is shown in Section \ref{sec:proofs-SB}, for any linear operator $T$ and $1 \leq k < {\rm rank}(T)=m$,
\begin{equation} \label{eq:sb-gaussian}
\left( \E \|T G\|_{2}^{-k} \right)^{-\frac{1}{k}}  = a_{k}(T)  \left( \E \|G_m\|_{2}^{-k} \right)^{-\frac{1}{k}} \sim a_{k}(T)\sqrt{m},
\end{equation}
where $G_m$ is the standard Gaussian random vector in $\R^m$.

That, combined with Theorem \ref{main-compl-sb}, leads to the accurate small-ball behaviour of $TX$:
\begin{Corollary} \label{cor:SB-intro}
There is an absolute constant $c$ such that the following holds. If $X$ satisfies the {\it SBA} with constant ${\cal{L}}$ and $ T: \mathbb R^{n} \rightarrow \mathbb R^{m}$ then for $ 1\leq k \leq m= {\rm rank}(T)$ and every $\eps>0$,
$$
\mathbb P \left ( \| T X \|_{2} \leq \varepsilon \sqrt{m} a_{k} (T)\right) \leq ( c {\cal{L}} \varepsilon)^{k}.
$$
\end{Corollary}
As it happens, one can control $ a_{k}(T) $ in terms of $ \|T\|_{S_{2}}$ as long as the operator $T$ does not have a trivial $q$-stable rank:

\begin{Theorem} \label{thm:small-ball-gaussian}
For every $q>2$ there are constants $c_q$ and $c_q^\prime$ that depend only on $q$, and absolute constants $c$ and $c^\prime$ such that the following holds.
Let $X$ be a random vector in $\mathbb R^{n}$ that satisfies the {\it SBA} with constant ${\cal{L}}$ and $ T: \mathbb R^{n} \rightarrow \mathbb R^{m}$  with $ m={\rm rank}(T)$.  For every $k \leq c_q {\rm srank}_{q} (T)$,
$$
\left( \E \|T G\|_{2}^{-k} \right)^{-\frac{1}{k}} \geq c\|T\|_{S_2};
$$
in particular, for every $\eps>0$
\begin{equation} \label{main-compl-sb-3}
\mathbb P \left ( \|T X\|_{2} \leq \frac{\eps}{2e\cal{L}} \| T \|_{S_{2}} \right) \leq (c^\prime \eps)^{c^\prime_q {\rm srank}_{q} (T) }.
\end{equation}
\end{Theorem}

The proofs of Theorem \ref{main-compl-sb} and Theorem \ref{thm:small-ball-gaussian} are presented in Section \ref{sec:proofs-SB}.

\subsection{Coordinate small ball}
As we noted previously, the fact that the Euclidean norm $\|TX\|_2$ is likely to be large gives limited information on the geometry of the random vector $TX$. Most notably, it says nothing on the crucial feature that leads to stable point separation---the number of large coordinates $TX$ has with respect to a fixed orthonormal basis $(u_i)_{i=1}^m$. The coordinate small-ball estimate we establish is based on the {\it wSBA}, and shows that indeed many of the coordinates  $(\inr{TX,u_i})_{i=1}^m$ are likely to be large. To see what sort of information on the coordinates $(\inr{TX,u_i})_{i=1}^m$ one can hope for, let us return to the Gaussian case (which, based on Theorem \ref{main-compl-sb}, is a likely candidate to be the `worst' random vector that satisfies the {\it wSBA}).

\begin{Example} \label{EX:CSB-independent}
Let $n=m$, set $X=G=(g_i)_{i=1}^m$ and consider the identity operator $T=Id:\R^m \to \R^m$. Given any orthonormal basis  $(u_i)_{i=1}^m$ it follows from rotation invariance and independence that
$$
\mathbb P \left( \left|\{i : |\inr{G,u_i}| \geq \eps\} \right| \leq c_1 m\right) = \mathbb P \left( \left|\{i : |g_i| \geq \eps\} \right| \leq c_1 m\right) \leq (c_2 \eps)^{c_3 m}
$$
for absolute constants $c_1,c_2$ and $c_3$.

Thus, the fact that $\|G\|_2$ is likely to be ~$\gtrsim \sqrt{m}$ is exhibited by a proportional number of the coordinates $(\inr{G,u_i})_{i=1}^m$ whose absolute values are larger than $\eps \|Id\|_{S_2}/\sqrt{m} = \eps$.
However, in general, obtaining a coordinate small-ball estimate is a nontrivial task even when $X$ has iid coordinates and $T$ is the identity operator. Indeed, let $n=m$ and set $X=(x_i)_{i=1}^m$ where the $x_i$'s are independent copies of a mean-zero random variable $x$. When $(u_i)_{i=1}^m$ is the standard basis, one has that $|\inr{TX,u_i}|=|x_i|$, and estimating
$$
\PROB \left( \sum_{i=1}^m \IND_{\{|x_i| \geq \eps\}} \leq \ell \right)
$$
is easy to do thanks to the independence of the $x_1,...,x_m$. But when $(u_i)_{i=1}^m$ is a different orthonormal basis then the coordinates of $(\inr{X,u_i})_{i=1}^m$ are likely to have strong dependencies and the wanted estimate is far from obvious.
\end{Example}

We present two coordinate small-ball estimates: Theorem \ref{th2-main-s2}, when the linear operator $T$ satisfies that $\|T^{\ast} u_i\|_2 =1$ for every $1 \leq i \leq m$, and Theorem \ref{th2-main-s2-gen} for more general operators $T$.

\begin{Theorem} \label{th2-main-s2}
There exists an absolute constant $c$ such that the following holds. Let $X$ satisfy the {\it wSBA} with constants $\theta$ and ${\cal L}$, set $(u_i)_{i=1}^m$ to be an orthonormal basis of $\R^m$ and consider $T:\R^n \to \R^m$ such that $\|T^{\ast} u_i\|_2=1$ for every $1 \leq i \leq m$. Let $q>2$ and set $ k_{q}= {\rm srank}_{q} (T) $. Then for $ s \in (0,1)$
\begin{equation} \label{th2-main-s2-1}
\PROB \left( \left| \left\{ i \leq m : | \inr{TX , u_{i}}| \geq\theta  \right\} \right|\leq (1-s) m \right) \leq 2  \left( \frac{2}{s}\right)^{ \frac{q}{ q-2}} \frac{m}{ k_{q}} \left( \frac{ c_{q}{\cal L}\theta}{ s}\right)^{ \frac{1}{2} (s/2)^{\frac{q}{q-2}} k_{q}  },
\end{equation}
where $c_{q} \leq c(q/(q-2))^{1/2}$.
\end{Theorem}

To put Theorem \ref{th2-main-s2} is some perspective, let us return to Example \ref{EX:CSB-independent}. Consider the case where $n=m$ and $T=Id$. Thus, $k_\infty=m$ and $\|T\|_{S_2}/\sqrt{m}=1$.  If $(u_i)_{i=1}^m$ is an arbitrary orthonormal basis and $s=0.01$, then by Theorem \ref{th2-main-s2} one has that with probability at least $1-2\left(c_1{\cal L}\eps\right)^{c_2 m}$,
\begin{equation} \label{eq:example-intro-1}
\left|\left\{i : |\inr{X,u_i}| \geq \eps\right\} \right| \geq 0.99m.
\end{equation}
Recall that by \cite{RV}, if $x$ is a random variable that has a density that is bounded by ${\cal L}$ and $X=(x_i)_{i=1}^m$ has iid coordinates distributed according to $x$, then $X$ satisfies the {\it SBA} with constant $c{\cal L}$; in particular, $X$ satisfies the {\it wSBA} with constants $\theta$ and $c{\cal L}$ for any $\theta>0$. Thus, with probability at least $1-2(c_1{\cal L}\eps)^{c_2m}$,
$$
\left| \left\{ i : |\inr{X,u_i} | \geq \eps \right\} \right| \geq 0.99m,
$$
which means that $X$ exhibits the same coordinate small-ball behaviour with respect to an \emph{arbitrary basis} as it would with respect to the standard basis; moreover, that behaviour is at least as good as that of the standard Gaussian vector.

\begin{Remark}
The one place in which Theorem \ref{th2-main-s2} is potentially loose is the factor $m/k_{q}$. It has an impact only in situations where the operator $T$ is, in some sense, trivial---when the $q$-stable rank of $T$ is smaller than $c_q \log m$.
\end{Remark}

The main application of Theorem \ref{th2-main-s2} is Theorem \ref{thm:circulant-small-ball}, showing that random sub-sampled convolutions exhibit stable point separation. In addition to that, a further application of Theorem \ref{th2-main-s2} is an $\ell_p$ small-ball estimate.
\begin{Theorem} \label{thm:l-p}
There exists absolute constants $c_1$ and $c_2$ such that the following holds. Let $X$ be a random vector in $\R^{n}$ that satisfies the {\it SBA} with constant ${\cal{L}}$. Set$ a \in \R^{n}$ and let $k=(c_1\|a\|_p/\|a\|_\infty)^{p}$. Then for any $0<\eps<1$,
\begin{equation}\label{eq:thm-lp-sb}
\PROB \left(\left\|\sum_{i=1}^m a_i x_i e_i\right\|_p \leq \eps \|a\|_p \right) \leq (c_2 \eps {\cal L})^{k}.
\end{equation}
\end{Theorem}

To see that \eqref{eq:thm-lp-sb} is truly a small-ball estimate with respect to the $\ell_p$ norm, observe that by the {\it SBA}
$$
\E \left\|\sum_{i=1}^m a_i x_i e_i\right\|_p^p = \sum_{i=1}^m |a_i|^p \cdot \E |x_i|^p \geq (c{\cal L})^p \|a\|_p^p;
$$
therefore, under a suitable moment assumption, $\|a\|_p \sim \E \left\|\sum_{i=1}^m a_i x_i e_i\right\|_p$.

\vskip0.4cm
There is no obvious way of obtaining an upper bound on \eqref{eq:thm-lp-sb}. If $X$ has iid coordinates and satisfies the {\it SBA} with constant ${\cal L}=1$, one may invoke \cite{PP-IMA}, where it is shown that for any semi-norm $\| \ \|$ and any $u>0$,
\begin{equation} \label{eq:PP-IMA}
\PROB(\|X\| \leq u) \leq \PROB(\|Y\| \leq u),
\end{equation}
and $Y$ is the uniform measure on $[-\frac{1}{2},\frac{1}{2}]^m$. However, similar comparison results of this kind for a general random vector $X$---whose coordinates need not be independent---are not known.

The proof of Theorem \ref{thm:l-p} is presented in Section \ref{sec:ell-p}.

\vskip0.4cm

We end the introduction with some notation. Throughout, $c, c_1, c^\prime$, etc., denote absolute constants. Their value may change from line to line. $c_q$ and $c(q)$ denote constants that depend on the parameter $q$; $a \lesssim b$ means that there is an absolute constant $c$ such that $a \leq cb$; and $a \lesssim_q b$ implies that $c$ depends on the parameter $q$. The corresponding two-sided estimates are denoted by $a \sim b$ and $a \sim_q b$ respectively.

For a subspace $F \subset \R^n$ let $P_F$ be the orthogonal projection onto $F$; $(e_i)_{i=1}^n$ is the standard basis of $\R^n$ and $P_k$ is the orthogonal projection onto ${\rm span}(e_1,...,e_k)$.
The standard Gaussian random vector in $\R^n$ is denoted by $G$, while $G_m$ is the standard Gaussian random vector in $\R^m$. Finally, if $f_X$ is the density of a random vector $X$, the density of $P_FX$ is denoted by $f_{P_FX}$.

\section{Proofs: small ball estimates} \label{sec:proofs-SB}
The starting point of the proof of Theorem \ref{main-compl-sb} is the following equality (see \cite{P-sb}, Proposition 4.6):

\begin{Proposition}  \label{P-sb}
 For every random vector $W$ in $ \mathbb R^{m}$ with bounded density and $ 1\leq k \leq m-1$,
 \begin{equation} \label{P-sb-1}
  \frac{ \left( \mathbb E \| W\|_{2}^{-k} \right)^{-\frac{1}{k}} }{  \left( \mathbb E \| G_m\|_{2}^{-k} \right)^{-\frac{1}{k}} } = \frac{1}{ \sqrt{2\pi}} \left( \int_{{\cal G}_{m,k}} f_{P_{F}W}(0)  d F \right)^{-\frac{1}{ k}}  ,
 \end{equation}
 with integration taking place with respect to the Haar measure on the Grassman manifold ${\cal G}_{m,k}$.
 \end{Proposition}

Proposition \ref{P-sb} indicates the path the proof of Theorem \ref{main-compl-sb} follows: one obtains suitable lower bounds on the $L_\infty$ norms the densities of typical projections of $TX$. This requires two straightforward volumetric observations that also explain the role of the Gaussian parameters $a_k(T)$.

\begin{Lemma} \label{SVD-lem}
Let $X$ be a random vector with a density. Consider $S:\R^n \to \R^k$ for $k\leq n$ and with ${\rm rank}(S)= k$, and let $UDP_kV$ be the singular value decomposition of $S$. Then, for any compact subset $K \subset \R^k$
\begin{equation} \label{SVD-lem-1}
\mathbb P \left( S X \in K \right) = \mathbb P \left( P_{E} X \in V^{\ast} D^{-1} U^{\ast} K \right),
\end{equation}
where $E=V^{\ast} (\R^k)$.

Moreover,
\begin{equation} \label{SVD-lem-2}
\mathbb P \left( S X \in K \right) \leq {\rm det } (D^{-1}) {\rm vol}(K)  \| f_{P_{E} X} \|_{L_\infty}  = \frac{\rm vol(K)}{\sqrt{ {\rm det} ( SS^{\ast} )}} \| f_{P_{E} X} \|_{L_\infty}.
\end{equation}
\end{Lemma}

\proof Since $S = UDP_k V$ and $P_k V = V P_E$, it follows that
\begin{align*}
\mathbb P \left( S X \in K \right)  = & \mathbb P \left( U D P_{k} V X \in K \right) = \mathbb P\left( P_{k} V X \in D^{-1} U^{\ast} K \right)
\\
= & \mathbb P \left( VP_{E} X \in D^{-1} U^{\ast} K \right) = \mathbb P \left( P_{E} X \in V^{\ast} D^{-1} U^{\ast} K \right)=(*);
\end{align*}
and by a volumetric estimate,
$$
(*) = \int_{V^{\ast}D^{-1} U^{\ast} K } f_{P_{E} X} ( x) d x \leq {\rm vol}(V^{\ast}D^{-1} U^{\ast} K) \| f_{P_{E} X} \|_{\infty} =  {\rm det } (D^{-1}) {\rm vol}(K) \| f_{P_{E} X} \|_{\infty}.
$$
\endproof

The second observation yields an estimate on the $L_\infty$ norm of the density of a projection of the random vector $TX$.

\begin{Lemma} \label{lem-bound-infty}
Let $X$ be a random vector, set $ 1\leq k \leq m-1\leq n-1$ and assume that for every $ E\in {\cal G}_{n,k}$,
\begin{equation} \label{as-1}
\| f_{P_{E} X}\|_{L_\infty} \leq {\cal{L}}^{k}.
\end{equation}
Then, for every $F \in {\cal G}_{m,k}$ and $T:\R^n \to \R^m$,
\begin{equation} \label{lem-gen-1}
\| f_{{P_F}T  X}\|_{L_\infty} \leq \frac{{\cal{L}}^{k} }{({\rm det} [( P_{F} T )(  P_{F} T )^{\ast}])^{\frac{1}{2}}}.
\end{equation}
\end{Lemma}

\proof Fix $F \in {\cal G}_{m,k}$ and observe that for every compact set $ K \subset F$,
$$
\frac{1}{{\rm vol}(K)} \int_{K} f_{P_{F}T X } ( x) d x  = \frac{1}{{\rm vol}(K)} \mathbb P \left( P_F T X \in K \right).
$$
By \eqref{SVD-lem-2} and the uniform estimate on $\| f_{P_{E} X}\|_{L_\infty}$ it follows that
$$
\mathbb P \left( P_F T X \in K \right) \leq {\rm vol}(K) \cdot \max_{E\in {\cal G}_{n,k}}\frac{\| f_{P_{E} X}\|_{L_\infty}}{ ({\rm det} [( P_F T )(  P_F T )^{\ast}])^\frac{1}{2} } \leq {\rm vol}(K) \cdot \frac{  {\cal{L}}^{k}}{ ({\rm det} [( P_F T )(  P_F T )^{\ast}])^\frac{1}{2} }.
$$
Therefore,
$$
\frac{1}{{\rm vol}(K)} \int_{K} f_{P_{F}T X } ( x) d x \leq \frac{  {\cal{L}}^{k}}{ ({\rm det} [( P_F T )(  P_F T )^{\ast}])^\frac{1}{2} },
$$
and since the R.H.S. is independent of $K$ the claim follows.

\endproof

\noindent {\bf Proof of Theorem \ref{main-compl-sb}.}
Recall that $G$ is the standard Gaussian random vector in $\R^n$ and $G_m$ is the standard Gaussian random vector in $\R^m$. The proof follows by invoking Proposition \ref{P-sb} twice: it is used to compare negative moments of $TG$ and $G_m$, and then to compare negative moments of $G_m$ and $TX$.

Let $F \in {\cal G}_{m,k}$. Since $P_FTG$ is also a centred Gaussian vector, it standard to verify that
\begin{equation}
\label{Gauss-st}
f_{P_{F} T G}^{\frac{1}{k}} (0) = \frac{1}{ \sqrt{2\pi}}  \left(\frac{ 1}{ {\rm det}[(P_{F} T) ( P_{F} T)^{\ast}]}\right)^{\frac{1}{2k}}.
\end{equation}
Hence, by \eqref{P-sb-1} and the definition of $a_k(T)$,
 \begin{align} \label{eq:ratio-gaussian}
 \frac{ \left( \E \| TG\|_{2}^{-k} \right)^{-\frac{1}{k}} }{  \left( \E \| G_m\|_{2}^{-k} \right)^{-\frac{1}{k}} } = & \frac{1}{ \sqrt{2\pi}} \left( \int_{{\cal G}_{m,k}} f_{P_{F}TG}(0)  d F \right)^{-\frac{1}{ k}} =  \left( \int_{{\cal G}_{m,k}} {\rm det}^{-\frac{1}{2}}[(P_{F} T) ( P_{F} T)^{\ast}]  d F \right)^{-\frac{1}{ k}}  \nonumber
 \\
 = & a_{k} ( T) .
 \end{align}
 On the other hand,
 $$
 \frac{ \left( \E \| TX\|_{2}^{-k} \right)^{-\frac{1}{k}} }{  \left( \E \|G_m\|_{2}^{-k} \right)^{-\frac{1}{k}} } = \frac{1}{ \sqrt{2\pi}} \left( \int_{{\cal G}_{m,k}} f_{P_{F} TX}(0)  d F \right)^{-\frac{1}{ k}};
 $$
by Lemma \ref{lem-bound-infty}, for every $F \in {\cal G}_{n,k}$,
$$
f_{P_{F} TX}(0) \leq \| f_{{P_F}TX}\|_{L_\infty} \leq \frac{{\cal{L}}^{k}}{ ({\rm det} [( P_{F} T )(  P_{F} T )^{\ast}])^\frac{1}{2}},
$$
implying that

\begin{equation*}
\frac{1}{ \sqrt{2\pi}} \left( \int_{{\cal G}_{m,k}} f_{P_{F} TX}(0)  d F \right)^{-\frac{1}{ k}} \geq \frac{1}{ \sqrt{2\pi}{\cal L}} \left( \int_{{\cal G}_{m,k}} {\rm det}^{-\frac{1}{2}} [( P_{F} T )(  P_{F} T )^{\ast}] dF \right)^{-\frac{1}{k}}
= \frac{a_k(T)}{\sqrt{2 \pi} {\cal L}}.
\end{equation*}
 Therefore,
 $$
 \frac{ \left( \E \| TX\|_{2}^{-k} \right)^{-\frac{1}{k}} }{  \left( \E \| TG/(\sqrt{2\pi}{\cal{L}})\|_{2}^{-k} \right)^{-\frac{1}{k}} }  \geq 1,
 $$
 as claimed.

\endproof

Note that \eqref{P-sb-1} and \eqref{Gauss-st} imply that $ \left( \mathbb E \| T G \|_{2}^{-k} \right)^{-\frac{1}{k}} = a_{k} (T) \left( \mathbb E \| G_{m} \|_{2}^{-k} \right)^{-\frac{1}{k}} $ as claimed in \eqref{eq:sb-gaussian}.

\subsection{Proof of Theorem \ref{thm:small-ball-gaussian}}
Thanks to Theorem \ref{main-compl-sb}, it suffices to obtain a suitable lower bound on
$(\E\|TG\|_2^{-k})^{-\frac{1}{k}}$ for $k \lesssim {\rm srank}_q(T)$ and $q>2$.

\begin{Lemma} \label{lemma:s-rank-HS}
There exists an absolute constant $c$ for which the following holds. Let $0<\theta<1$ and $q>2$, and set
$$
m=(c\theta)^{2q/(q-2)}{\rm srank}_q(T).
$$
If $(s_i)_{i=1}^r$ are the non-zero singular values of $T$ arranged in a non-increasing order and
$$
\tilde{s_i} = \min\{s_i,\|T\|_{S_2}/\sqrt{m}\}
$$
then
$$
\sum_{i=1}^r \tilde{s}_i^2 \geq (1-\theta^2)\|T\|_{S_2}^2.
$$
\end{Lemma}

\proof Let $0<\theta<1$ and observe that for every $1 \leq i \leq m$, $s_i \leq \|T\|_{S_q}/i^{1/q}$. Therefore,
$$
\sum_{i \leq m} s_i^2 \leq \|T\|_{S_q}^2 \sum_{i \leq m} \frac{1}{i^{2/q}} \leq \frac{cq}{q-2} \|T\|_{S_q}^2 m^{1-2/q} \leq \theta^2 \|T\|_{S_2}^2
$$
provided that
$$
m \leq (c\theta)^{2q/(q-2)} {\rm srank}_q(T).
$$
At the same time, $s_{m+1} \leq \|T\|_{S_2}/\sqrt{m+1}$, implying that
$$
|\{i : s_i \geq \|T\|_{S_2}/\sqrt{m} \}| \leq m.
$$
Thus,
$$
\sum_{i=1}^r \min \left\{ \|T\|_{S_2}/\sqrt{m}, s_i\right\}^2 \geq \sum_{i=m+1}^r s_i^2 \geq (1-\theta^2)\|T\|_{S_2},
$$
and the claim follows. 
\endproof

The proof of Theorem \ref{thm:small-ball-gaussian} is based on the following outcome of the so-called ``B-Theorem" (see \cite{CFM} for the proof of the ``B-Theorem") and requires some additional notation.

\vskip0.4cm

For $a \in \R^n$ let $G_a = \inr{G,a}$, and for $A \subset \R^n$ set
$$
d_{\ast}(A) = \left(\frac{\E \sup_{a \in A} G_a}{\sup_{a \in A} (\E G_a^2)^{1/2}}\right)^2.
$$

\begin{Theorem} \label{thm:B-thm} \cite{KV,LO}.
There are absolute constants $c_1$ and $c_2$ such that for any $A \subset \R^n$ and any $0<s<1$,
$$
\mathbb P \left( \sup_{a \in A} G_a \leq s \E \sup_{a \in A} G_a \right) \leq (c_1s)^{c_2 d_{\ast}(A)}.
$$
\end{Theorem}

\noindent{\bf Proof of Theorem \ref{thm:small-ball-gaussian}.}
Let $(s_i)_{i=1}^r$ be the non-zero singular values of $T$ and set $(\tilde{s}_i)_{i=1}^r$ to be as in the proof of Lemma \ref{lemma:s-rank-HS}. Using the notation of the lemma, let $\theta^2=3/4$. Note that if $D$ is a diagonal operator that satisfies  $d_{ii}=s_i$ for $i \leq r$ and $0$ otherwise, and $\tilde{D}$ is a diagonal operator whose non-zero diagonal entries are $d_{ii}=\tilde{s_i}$ for $i \leq r$, then
$$
\tilde{D} B_2^n \subset D B_2^n, \ \ \ \|\tilde{D}\|_{S_\infty} \leq \frac{\|T\|_{S_2}}{\sqrt{m}}, \ \ {\rm and} \ \  \|\tilde{D}\|_{S_2} \geq \frac{\|T\|_{S_2}}{2}.
$$

By rotation invariance, for every $k$, $\E \|TG\|_{2}^{-k} = \E\|DG\|_2^{-k}$, and for every $x \in \R^n$, $\|\tilde{D}x\|_2 \leq \|Dx\|_2$. Hence,
$$
(E \|TG\|_2^{-k})^{-\frac{1}{k}} = (\E\|DG\|_2^{-k})^{-\frac{1}{k}} \geq (\E\|\tilde{D}G\|_2^{-k})^{-\frac{1}{k}}.
$$

Let $A =\tilde{D} B_2^n$ and observe that for $t \in \R^n$,
$$
\sup_{a \in A} \inr{a,t} = \sup_{x \in B_2^n} \inr{x,\tilde{D}t} = \|\tilde{D}t\|_2;
$$
therefore,
$$
\E \sup_{a \in A} \inr{a,G} = \E \|\tilde{D}G\|_2 \geq \|\tilde{D}\|_{S_2} \geq  \frac{\|T\|_{S_2}}{2}
$$
and
$$
\sup_{a \in A} \E \inr{a,G}^2 \leq \max_i d_{ii} \leq \frac{\|T\|_{S_2}}{\sqrt{m}}.
$$
Finally, by Theorem \ref{thm:B-thm}, for every $0<u<1$,
$$
\mathbb P ( \|\tilde{D}G\|_2 \leq c_1u\|T\|_{S_2} ) \leq (u/2)^{c_2m},
$$
where $c_1$ and $c_2$ are suitable absolute constants. A straightforward tail integration argument shows that for $k \leq c_3m$
$$
(\E\|\tilde{D}G\|_2^{-k})^{-\frac{1}{k}} \geq c_4\|T\|_{S_2},
$$
as required.

\section{Proofs: Coordinate small-ball estimates}
Let us turn to the proof of Theorem \ref{th2-main-s2}. Recall that $X$ is an $n$-dimensional random vector that satisfies the {\it wSBA} with constants $\theta$ and ${\cal L}$, let $ (u_i)_{i=1}^m$ be an orthonormal basis of $\R^{m}$ and set $ T:\R^{n} \to \R^{m}$ to be a linear operator which satisfy that for $1 \leq i \leq m$
\begin{equation} \label{T-Assum}
\| T^{\ast} u_{i} \|_{2} = 1.
\end{equation}

The key component of the proof of Theorem \ref{th2-main-s2} is a decomposition lemma. To formulate it, let $\sigma \subset \{1,...,m\}$ and denote by $P_\sigma : \R^m \to \R^\sigma$ the orthogonal projection onto ${\rm span}(u_i)_{i \in \sigma}$. Thus, $P_{\sigma}^{\ast}: \R^{\sigma} \to \R^m$ is the formal identity operator with respect to the basis $(u_i)_{i=1}^m$.

\begin{Lemma} \label{NY-lem}
Let $q>2$ and set $c_q \sim (q/(q-2))^{1/2}$. Assume that for every $1 \leq i \leq m$, $\|T^{\ast}u_i\|_2=1$ and set $k_{q}= {\rm srank}_{q}(T)$. Then for any $\lambda \in (0,1)$ there are disjoint subsets $ \sigma_{1},..., \sigma_{\ell} \subset \{1,...,m\}$ such that
\begin{description}
\item{$\bullet$} For $1 \leq j \leq \ell$, $| \sigma_{j} | \geq  \lambda^{\frac{q}{q-2}} k_{q}/2$ and $\sum_{j=1}^{\ell}|\sigma_{j}| \geq (1-\lambda) m  $; and
\item{$\bullet$} $\| (T^{\ast}P^{\ast}_{\sigma_j})^{-1}\|_{S_{\infty}} \leq  c_q$.
\end{description}
\end{Lemma}

The proof of Lemma \ref{NY-lem} is based on the idea of restricted invertibility. The version used here is Theorem 8 from  \cite{NY}:

\begin{Theorem} \label{NY} 
For $q>2$ set $c_q \sim (q/(q-2))^{1/2}$. If $ A:\R^{m} \to \R^{n}$ is a linear operator then there exists $ \sigma \subset \{1,...,m\}$ of cardinality at least $| \sigma| \geq {\rm srank}_{q}(A)/2$ such that the map $ (AP^{\ast}_{\sigma})^{-1}$ is well defined and
\begin{equation} \label{NY-1}
\| (AP^{\ast}_{\sigma})^{-1} \|_{S_{\infty}} \leq c_q \frac{ \sqrt{m}}{ \| A\|_{S_{2}} }.
\end{equation}
\end{Theorem}

\noindent {\bf Proof of Lemma \ref{NY-lem}.} The construction of the subsets $(\sigma_j)_{j=1}^\ell$ is performed inductively. First, let $k_q = {\rm srank}_q(T)$ and apply Theorem \ref{NY} to $A=T^{\ast}$. Thus, noting that $\|T^*\|_{S_2} = \sqrt{m}$, there is $ \sigma_1 \subset \{1,...,m\}$ such that
\begin{equation} \label{NY-2-1}
 |\sigma_{1}|\geq \frac{k_q}{2} \ \ \  {\rm and} \  \  \ \| ( T^{\ast}P^{\ast}_{\sigma_1} )^{-1}\|_{S_{\infty}} \lesssim \left(\frac{q}{q-2}\right)^{1/2}.
\end{equation}

If $|\sigma_{1}|\geq (1-\lambda) m$ the lemma is proved. Otherwise, let $m_1=m-|\sigma_1|$ and set $T_1= P_{\sigma_{1}^{c} } T  : \R^{n} \to \R^{m_{1}}$. Since $ \|T_{1}^{\ast} u_{i} \|_{2} =1$ for all $i\in \{1,...,m\} \setminus \sigma_{1}$, it is evident that $\|T_1\|_{S_2}^2 \geq \lambda m$, and because $P_{\sigma_1^c}$ is a contraction one has that $\|T_1\|_{S_q} = \|P_{\sigma_1^c} T\|_{S_q} \leq \|T\|_{S_q}$. Set $k_q^{(1)}= {\rm srank_{q}} ( T_{1})$; thus,
\begin{equation*}
k_{q}^{(1)} = \left( \frac{ \| T_{1} \|_{S_{2}}}{ \|T_{1}\|_{S_{q}}} \right)^{\frac{2q}{ q-2}} \geq \lambda^{\frac{q}{q-2}}\left( \frac{ \| T \|_{S_{2}}}{ \|T\|_{S_{q}}} \right)^{\frac{2q}{ q-2}} =  \lambda^{\frac{q}{q-2}} k_{q}.
\end{equation*}
Invoking Theorem \ref{NY} again, this time for $A=T_1^{\ast}$, there is $ \sigma_{2} \subset \{1,...,m\} \setminus \sigma_{1}$, such that
\begin{equation*}
 |\sigma_2| \geq \frac{k_{q}^{(1)}}{2} \geq \frac{1}{2}\lambda^{\frac{q}{q-2}} k_{q} \ \ \ {\rm and} \ \ \ \| ( T^{\ast} P^{\ast}_{\sigma})^{-1}\|_{S_{\infty}} \lesssim \left(\frac{q}{q-2}\right)^{1/2}.
 \end{equation*}

Again, if $|\sigma_{1} | + |\sigma_{2}| \geq (1-\lambda) m$ the lemma is proved, and if not one may continue in the same way, constructing operators $T_{j}$ and sets $ \sigma_{j}$ inductively until  $\sum_{j=1}^{\ell} |\sigma_{j}|\geq(1-\lambda) m$.
\endproof

\vskip0.4cm

The fact that $\| (T^{\ast}P^{\ast}_{\sigma})^{-1} \|_{S_{\infty}} \leq \gamma$ implies that the ellipsoid $P_\sigma T(B_2^n)$ contains the Euclidean ball $\gamma^{-1} B_2^\sigma$, which leads to a small-ball estimate.

\begin{Lemma} \label{lemma:decom-to-CSM}
There is an absolute constant $c$ for which the following holds.
Let $\sigma \subset \{1,...,m\}$ such that $\| (T^{\ast}P^{\ast}_{\sigma} )^{-1}\|_{S_{\infty}} \leq \gamma$, and set $ \gamma_{0}= \max\{ 1, \gamma\}$.
If $X$ satisfies the {\it wSBA} with constants $\theta$ and ${\cal L}$, then for any $\tau \subset \sigma$,
$$
\mathbb P (\|P_\tau T\|_2 \leq\theta \sqrt{|\tau|} ) \leq (c \gamma_{0} \theta {\cal L})^{|\tau|}.
$$
\end{Lemma}

An observation one needs for the proof of Lemma \ref{lemma:decom-to-CSM} is a monotonicity property for the {\it wSBA}. Its proof can be found in Proposition 2.1 in \cite{RV} and is based on a simple covering argument.
\begin{Lemma} \label{lemma:wSBA-M}
Let $X$ satisfy the {\it wSBA} with constants $\theta$ and ${\cal L}$. Then for every $M >1$, $X$ also satisfies the {\it wSBA} with constants $M \theta$ and $3 {\cal L}$.
\end{Lemma}

\noindent {\bf Proof of Lemma \ref{lemma:decom-to-CSM}.} Let $ r= |\tau|$ and note that $B=P_{\tau} P_{\sigma} T : \R^n \to \R^r$ is a linear operator of rank $r$. As noted previously, there are $|\sigma|$ non-zero singular values of $P_\sigma T$, all of which are at least $\gamma^{-1}$; and since $\tau \subset \sigma$, it follows that $ s_{i}(B) \geq \gamma^{-1}$ for $ 1\leq i \leq r$. By the singular value decomposition theorem there are $ U \in {\cal{O}}_{r}$, $ V\in {\cal{O}}_{n}$ and a diagonal matrix $ D= {\rm diag}( s_{1}(B), \cdots , s_{r}(B))$ such that
$B= U D P_{r} V$, where, as always, $ P_{r}$ denotes the orthogonal projection onto $\{e_1,...,e_r\}$. Setting $ F= V^{\ast}(\R^{r}) $ one has that
$P_{r} V = VP_{F}$.

Since $X$ satisfies the {\it wSBA} with constants $\theta$ and ${\cal L}$, and all the entries in the diagonal of $D$ are at least $\gamma^{-1}$, invoking Lemma \ref{lemma:wSBA-M} it is evident that
\begin{align*}
\mathbb P \left( \| P_{\tau} T X \|_{2}\leq \theta \sqrt{ | \tau|} \right) = & \mathbb P \left( \| UDVP_{F} X \|_{2}\leq\theta \sqrt{ | \tau|} \right) \leq
\mathbb P \left( \| P_{F} X \|_{2}\leq \gamma \theta \sqrt{ | \tau|} \right)
\\
\leq  & \left( \theta \gamma_{0} \cal L   \right)^{ | \tau|}  ,
\end{align*}
and the claim follows.
\endproof

\vskip0.4cm

\noindent {\it Proof of Theorem \ref{th2-main-s2}}. Let $(\sigma_j)_{j=1}^\ell$ be the collection of subsets as in Lemma \ref{NY-lem} and set
$$
\sigma= \bigcup_{j\leq \ell }~ \sigma_j.
$$
In particular, for $1 \leq j \leq \ell$,
$$
|\sigma_j| \geq \lambda^{\frac{q}{q-2}}k_q/2.
$$
Recall that $(u_i)_{i=1}^m$ is an orthonormal basis of $\R^m$ and for $\tau \subset \{1,...,m\}$ set
$$
Q_\tau = \{x \in \R^m : \max_{i \in \tau} |\inr{x,u_i}| \leq 1\}.
$$
For $1 \leq i \leq m$ consider the random variables
\begin{equation} \label{eta-def}
 \eta_i= \IND_{\{ z: | \inr{T z, u_i} | \geq \theta  \}} (X) \ \ \ \ \ {\rm and} \ \ \ \ \ \zeta_i= \IND_{\{ z: | \inr{T z, u_i} | < \theta  \}} (X).
\end{equation}
Thus, for $0<\alpha<1$ and $1 \leq j \leq \ell$,
\begin{align*}
& \PROB \left( \sum_{i\in \sigma_j } \zeta_{i} \geq \alpha | \sigma_j | \right) \leq \sum_{\tau \subset \sigma_j, \ | \tau|= \alpha | \sigma_j |} \PROB \left( \bigcap_{i\in \sigma_j} \{  | \inr{T X, u_i}| < \theta  \} \right)
\\
\leq & \left( \frac{e}{ \alpha } \right)^{ \alpha | \sigma_{j}|} \max_{\tau \subset \sigma_j, \ | \tau|= \alpha | \sigma_j| } \PROB \left( P_{\tau} T X \in \theta Q_{\tau}  \right) \leq  \left( \frac{e}{ \alpha } \right)^{ \alpha | \sigma_{j}|} \max_{\tau \subset \sigma_{j}, \ | \tau|= \alpha | \sigma_j| } \PROB \left( \| P_{\tau} T X \|_{2} <  \theta \sqrt{| \tau|}  \right),
\end{align*}
where the last inequality holds because $Q_\tau \subset \sqrt{|\tau|}B_2^\tau= \sqrt{|\tau|} B_{2}^{m} \subset \mathbb R^{\tau} $. Since $X$ satisfies the {\it wSBA} with constants $\theta$ and ${\cal L}$ it follows from Lemma \ref{lemma:decom-to-CSM} that
$$
\max_{\tau \subset \sigma_{j}, \ | \tau|= \alpha | \sigma_j| } \mathbb P \left( \| P_{\tau} T X \|_{2} <  \theta \sqrt{| \tau|} \right) \leq \left( c_q {\cal L} \theta\right)^{\alpha |\sigma_{j}|}
$$
and $c_q \sim (q/(q-2))^{1/2}$; therefore,
\begin{equation} \label{th2-main-0-1}
 \mathbb P \left( \sum_{i\in \sigma_j } \eta_{i} \leq ( 1-\alpha) | \sigma_j | \right) \leq  \left( \frac{ ec_q {\cal L}\theta}{\alpha}\right)^{\alpha|\sigma_{j}|}.
\end{equation}
Now set $0<s<1$, let $\lambda=s/2$ and recall that $ |\sigma| \geq (1-\lambda) m $.  Set $1-\alpha = (1-s)/(1-\lambda)$ and note that $\alpha \geq s/2$. With this choice of $\lambda$, the union bound and  \eqref{th2-main-0-1},
\begin{align*}
& \PROB \left( \sum_{i=1}^{m} \eta_{i} \leq (1-s) m\right)  = \PROB \left(\sum_{i \in \sigma} \eta_i \leq (1-s) \cdot \frac{|\sigma|}{1-\lambda}\right)
\\
\leq & \PROB \left( \sum_{j\leq \ell } \sum_{i\in \sigma_{\ell}}  \eta_{i} \leq (1-\alpha)  \sum_{j\leq \ell } | \sigma_j | \right) \leq \sum_{j\leq \ell} \PROB \left(  \sum_{i\in \sigma_j}  \eta_{i} \leq (1-\alpha) | \sigma_j | \right)
\\
\leq  & \sum_{j\leq \ell} \PROB \left(  \sum_{i\in \sigma_j}  \eta_{i} \leq \left(1-\frac{s}{2}\right)  | \sigma_j | \right) \leq \sum_{j\leq \ell}  \left( \frac{ec_q {\cal L}\theta}{s/2} \right)^{|\sigma_{j}|/(s/2)}=(*).
\end{align*}
Finally, since
$$
|\sigma_j| \geq c_1 \lambda^{\frac{q}{q-2}}k_q \sim s^{q/(q-2)}k_q,
$$
it is evident that $\ell \lesssim s^{-q/(q/2)}m/k_q$ and the claim follows.
\endproof

\subsection{Coordinate small-ball for general operators}
The assumption that $\| T^{\ast} u_{i} \|_{2}=1$ for every $1 \leq i \leq m$
is not essential and can be replaced by a considerably weaker condition. If instead one assumes that there are constants $ \delta_{1}>0$ and $\delta_{2} \geq 1$ such that
 \begin{equation} \label{ass-NY-3}
\left( \frac{1}{ m} \sum_{i=1}^{m} \| T^{\ast} u_{i} \|_{2}^{2+\delta_{1}} \right)^{\frac{1}{ 2+ \delta_{1}}} \leq \delta_{2} \frac{\|T \|_{S_{2}}}{ \sqrt{m}},
\end{equation}
then the following version of Theorem \ref{th2-main-s2} can be established:

\begin{Theorem} \label{th2-main-s2-gen}
Let $X$ satisfy the {\it SBA} with constant ${\cal L}$. Consider $T:\R^n \to \R^m$, an orthonormal basis $(u_i)_{i=1}^m$ of $\R^m$ such that \eqref{ass-NY-3} is satisfied, and for $q> 2$ set $ k_{q}= {\rm srank}_{q} (T) $. Then, for any $\eps\in (0,1)$ one has that
\begin{equation} \label{th2-main-s2-1-gen}
\mathbb P \left( \left| \left\{ i \leq m : | \inr{TX , u_{i}}| \geq\eps  \frac{ \|T \|_{S_{2}}}{ \sqrt{m}}\right\} \right| \leq c_{0} m 	\right) \leq c_1 \frac{m}{k_q} (c_2 {\cal L} \eps)^{c_3 k_q}
\end{equation}
where $c_0=c_0(\delta_1,\delta_2)$, 
$$
c_1=  5^{ \frac{q}{ q-2}}, \ \ \ c_2 \sim  \left(\frac{q}{q-2}\right)^{1/2} \ \ \ {\rm and} \ \ c_3 = \left(\frac{1}{2\delta_2}\right)^{\frac{2+\delta_{1}}{ \delta_{1}}} \in (0,1).
$$
\end{Theorem}
\vskip0.4cm

Because the proof of Theorem \ref{th2-main-s2-gen} follows a similar path to that of Theorem \ref{th2-main-s2} we will only outline the necessary modifications.
\vskip0.4cm

\noindent{\bf Sketch of proof.} Note that \eqref{ass-NY-3} and the Paley-Zygmound inequality imply that
\begin{equation} \label{eq:in-proof-gen-CSB-1}
\left| \left\{ i \leq m : \| T^{\ast} u_{i} \|_{2} \geq \frac{ \|T \|_{S_{2}} }{2\sqrt{m}} \right\} \right| \geq c_{0}(\delta_1,\delta_2) m,
\end{equation}
and without loss of generality one may assume that $ c_{0} m $ is an integer. In particular, let $\sigma_{0}\subset \{1,...,m\}$ to be of cardinality $c_0 m$ and for every $i \in \sigma_0$,
$$
\|T^{\ast} u_i\|_2 \geq \frac{\|T\|_{S_2}}{2 \sqrt{m}}.
$$
We may assume that $\sigma_0=\{1,...,c_0m\}$ and let $ T_{0}= P_{\sigma_{0} } T$ where $P_{\sigma_0}$ is the orthogonal projection onto ${\rm span}(u_i)_{i \in \sigma_0}$.
Therefore,
$$
\|T \|_{S_{2}} \geq \| T_{0} \|_{S_{2}} \geq \frac{\sqrt{c_{0}}}{2} \| T \|_{S_{2}},
$$
and for any $ \sigma \subset \sigma_{0}$,
$$
\| P_{\sigma} T_{0} \|_{S_{2}} \geq \frac{1}{2}  \sqrt{\frac{ |\sigma|}{m}} \| T \|_{S_{2}} \geq \frac{ \sqrt{c_{0}} }{2}\sqrt{\frac{ |\sigma|}{|\sigma_{0}|}} \| T_{0} \|_{S_{2}}.
$$
Let $\delta = \sqrt{c_0}/2$ and set
$$
k_{q,0} = {\rm srank}_{q}(T_{0}) \geq \left( \frac{c_{0}}{ 4}\right)^{\frac{q}{q-2}} k_{q}.
$$
Following the argument used in the proof of Lemma \ref{NY-lem}, it is evident that there are disjoint subsets $ \sigma_{1},..., \sigma_{\ell} \subset \{1,...,c_0m\}$ such that
\begin{description}
\item{$\bullet$} For $1 \leq j \leq \ell$, $| \sigma_{j} | \geq (\delta^{2} \lambda)^{\frac{q}{q-2}} k_{q,0}/2$ and $\sum_{j=1}^{\ell}|\sigma_{j}| \geq (1-\lambda) |\sigma_0|=(1-\lambda)c_0m  $; and
\item{$\bullet$} $\| (T^{\ast}P^{\ast}_{\sigma_j})^{-1}\|_{S_{\infty}} \leq  \frac{c_0}{ \delta^{2}} $.
\end{description}
From here on the proof is identical to that of Theorem \ref{th2-main-s2} with the choice of $\lambda=1/2$; the details are omitted.
\endproof

\begin{remark}
A version of Theorem \ref{th2-main-s2-gen} holds true under the {\it wSBA} as well. We leave the details of the proof to the reader.
\end{remark}

\vskip0.4cm

Theorem \ref{th2-main-s2} and Theorem \ref{th2-main-s2-gen} imply that under mild assumptions on $T$, the coordinate small-ball estimate exhibits the standard small-ball one. Indeed, as an example, set $s=1/2$, let $k_4$ denote the $q$ stable rank for $q=4$ and observe that if $\|T^*e_i\|_2 =1$ for every $1 \leq i \leq m$ then $\|T\|_{S_2}=\sqrt{m}$. Hence,
$$
\left| \left\{ i : |\inr{TX,u_i}| \geq (\theta/\sqrt{2})\frac{\|T\|_{S_2}}{\sqrt{m}} \right\} \right| \geq \frac{m}{2}
$$
with probability at least
$$
1-(m/k_4) \cdot(c \theta {\cal L})^{c^\prime k_4}
$$
where $c$ and $c^\prime$ are absolute constants. Therefore, if $k_4 \gtrsim \log m$ and $\theta \lesssim 1/{\cal L}$,
$$
\mathbb P \left( \|TX\|_2 \leq \frac{\theta}{8} \|T\|_{S_2} \right) \leq (c^{\prime \prime} \theta {\cal L})^{c^{\prime} k_4/2},
$$
which recovers the small-ball estimate (and obviously similar bounds hold for any $q>2$ at the price of modified constants).

At the same time, the difference between the two estimates cannot be overstated: Theorem \ref{th2-main-s2} implies that for \emph{any} choice of a coordinate basis $(u_i)_{i=1}^m$, a typical realization of the vector $(\inr{TX,u_i})_{i=1}^m$ will have $\sim m$ large coordinates, which is a significantly stronger statement than the standard small-ball estimate. Indeed, there are many examples in which the coordinate structure dictated by the orthonormal basis is a feature of the problem, and a small-ball estimate is simply not good enough. The case of random sub-sampled convolutions, which we now turn to, is one such example.

\section{Proofs: applications}
Here we present the proofs of the applications that follow from Theorem \ref{th2-main-s2}, starting with point separation and stable point separation for random sub-sampled convolutions.

\subsection{Random sub-sampled convolutions and stable point separation} \label{sec:circulant}
Recall that $ \xi$ is an isotropic random vector in $\R^n$ and that $(\delta_i)_{i=1}^n$ are independent $ \{ 0,1 \}$-valued random variables with mean $ \delta$. If $I=\{i : \delta_i = 1\}$, the question of point separation is whether with high probability
\begin{equation} \label{eq:part-conv}
\frac{1}{\delta n} \sum_{i \in I} (a \circledast \xi)_i^2 \geq c_0 \|a\|_2^2
\end{equation}
for a suitable constant $c_0$ that is independent of $a$ and of $\delta$; and, as far as stable point separation is concerned, whether
\begin{equation} \label{eq:part-conv-stable}
\left|\left\{ i \in I : (a \circledast \xi)_i \geq c_1 \|a\|_2 \right\} \right| \geq c_2 \delta n.
\end{equation}

To analyze \eqref{eq:part-conv} and \eqref{eq:part-conv-stable}, let us first examine the convolution $a \circledast \xi$. Recall that ${\cal F}$ is the discrete Fourier matrix and observe that
$$
a \circledast \xi = {\cal F}^{-1} {\cal F} (a \circledast \xi) = {\cal F}^{-1} \left(({\cal F}a)_i \cdot ({\cal F} \xi)_i\right)_{i=1}^n = {\cal F}^{-1} D_{{\cal F}a} {\cal F} \xi,
$$
where $D_{{\cal F}a}$ is a diagonal matrix whose diagonal entries are $d_{ii}=({\cal F}a)_i$. Setting $U={\cal F}^{-1}/\sqrt{n}$, $W={\cal F}/\sqrt{n}$ and $O={\cal F}/\sqrt{n}$, it follows that
\begin{equation} \label{eq:circulant-via-orthogonal}
a \circledast \xi = \sqrt{n} UD_{Wa}O\xi \equiv \Gamma_a \xi,
\end{equation}
where $U,W$ and $O$ are Hadamard matrices.

Note that $\xi$ satisfies the {\it SBA} with constant ${\cal L}$; that $\hat{a}=\frac{{\cal F}a}{\sqrt{n}}=Wa$; and that $O$ and $U$ are orthogonal matrices. Hence, for every $q>2$,
$$
\|\Gamma_a\|_{S_q}= \sqrt{n}\|UD_{Wa}O\|_{S_q} = \sqrt{n} \|D_{Wa}\|_{S_q} = \sqrt{n} \|\hat{a}\|_{q}.
$$
Therefore, if $\|a\|_2=1$,
$$
{\rm srank}_{q}(\Gamma_a)=\left(\frac{\|\hat{a}\|_2}{\|\hat{a}\|_q}\right)^{\frac{2q}{q-2}}=
\left(\frac{1}{\|\hat{a}\|_q}\right)^{\frac{2q}{q-2}}.
$$

By Theorem \ref{th2-main-s2} for $(u_i)_{i=1}^n=(e_i)_{i=1}^n$, there is an event ${\cal A}$ of probability at least
$$
1-2  \left( c_1\right)^{ \frac{q}{ q-2}} \frac{n}{ k_{q}} \left( c_2(q){\cal L}\eps\right)^{ c_3^{\frac{q}{q-2}} k_{q}  }=(*)
$$
with respect to $\xi$, on which
$$
\left|\left\{ i : |\inr{\Gamma_a \xi ,e_i}| \geq \eps \right\} \right| \geq 0.99n.
$$
To control the probability estimate, let $W_1,...,W_n$ be the rows of the matrix $W$ and note that
$$
\max_{1 \leq i \leq n} \|W_i\|_{\infty}=\frac{1}{\sqrt{n}}.
$$
If $a$ is $s$-sparse with respect to the standard basis then
$$
\|\hat{a}\|_q^q = \sum_{i=1}^n |\inr{W_i,a}|^q \leq n \cdot \left(\frac{s}{n}\right)^{q/2}=\frac{s^{q/2}}{n^{(q-2)/2}},
$$
implying that
$$
(*) \geq 1- (c_4(q)\eps {\cal L})^{c_5(q)/\|\hat{a}\|_q^{2q/(q-2)}}
$$
as long as $s \leq c_6(q)n/\log n$ and $\eps \leq c_7(q)$.

Finally, for every realization of $\xi$ in the event ${\cal A}$, with probability at least $1-2\exp(-c_8 \delta n)$ with respect to $(\delta_i)_{i=1}^n$ one has that
$$
\left|\left\{ i \in I : |\inr{\Gamma_a \xi ,e_i}| \geq \eps \right\} \right| \geq 0.98\delta n,
$$
and in particular,
$$
\sum_{i \in I} |\inr{\Gamma_a \xi,e_i}|^2 \geq 0.98 \eps^2 \delta n.
$$
A Fubini argument completes the proof.
\endproof

\subsection{Small-ball estimates for the $\ell_p$-norm} \label{sec:ell-p}
Let $X=(x_i)_{i=1}^n$ satisfy the {\it SBA} with constant ${\cal L}$ and fix $a \in \R^n$. The goal here is to use Theorem \ref{th2-main-s2} and the information it provides on the distribution of the coordinates of a random vector $X$ to control the probability
$$
\PROB \left( \left\|\sum_{i=1}^n a_i x_i e_i\right\|_p \geq \eps \|a\|_p \right);
$$
here, as always, $(e_i)_{i=1}^n$ denotes the standard basis in $\R^n$ and $\| \ \|_p$ is the $\ell_p$ norm.

Without loss of generality assume that $a_1 \geq a_2 .... \geq 0$ and set
\begin{equation} \label{SB-1-sets-1}
I_{j} = \left\{ i : \frac{1}{ 2^{j+1}} < \frac{ a_i}{ a_1} \leq \frac{1}{2^{j}}\right\}.
\end{equation}
For every integer $\ell$ let
\begin{equation} \label{SB-1-sets-2}
\Lambda_{\ell} = \{ j : | I_j | = \ell\}
\end{equation}
and note that it is possible that some of the sets $ \Lambda_{\ell}$'s are empty. For every $\ell$ define
\begin{equation} \label{SB-1-sets-3}
j(\ell) = \min\Lambda_{\ell}
\end{equation}
and if $\Lambda_\ell$ is empty let $j(\ell)= 0$.

\vskip0.4cm

The idea behind this decomposition of $\{1,...,n\}$ to the union of the sets $I_j$ is that if $\Lambda_\ell \not = \emptyset$ then the contribution to $\|a\|_p$ that comes from $\bigcup_{j \in \Lambda_\ell} I_j$ is equivalent to the contribution of $I_{j(\ell)}$. Indeed, for any $I_j$,

 \begin{equation} \label{SB-1-sets-4}
\frac{1}{2^{p}}  | I_{j}|   \frac{a_{1}^{p}}{2^{jp}} < \sum_{i\in I_{j}} a_i^{p} \leq | I_{j}|   \frac{a_1^p}{2^{jp}}
\end{equation}
and by comparing the sum to an appropriate geometric progression, there are absolute constants $c_1$ and $c_2$ such that
$$
c_1^p \sum_{i \in I_{j(\ell)}} a_i^p \leq \sum_{j \in \Lambda_\ell} \sum_{i \in I_j} a_i^p \leq c_2^p \sum_{i \in I_{j(\ell)}} a_i^p.
$$
As a result, there are disjoint coordinate blocks, each one of different cardinality, such that

$$
\|a\|_p^p \sim \sum_{\{\ell \geq 1 : \Lambda_\ell \not = \emptyset\}} \sum_{i \in I_{j(\ell)}} a_i^p.
$$

Fix an index $\ell$ such that $\Lambda_\ell \not = \emptyset$ and consider $j=j(\ell)$. One has that
$$
\left\{i \in I_j : |x_i a_i| \leq  \eps \frac{a_1}{2^{j+1}}\right\} \subset  \left\{i \in I_j : |x_i| \leq  \eps \right\},
$$
and by Theorem \ref{th2-main-s2} for the orthogonal projection onto ${\rm span}(e_i: i \in I_j)$, denoted in what follows by $P_{I_j}$, there are absolute constants $c_3$ and $c_4$ such that
$$
\PROB\left(\left|\left\{i \in I_{j(\ell)} : |x_i| \leq  \eps \right\}\right| \leq \frac{\ell}{2} \right) \leq (c_3 {\cal L}\eps)^{c_4 \ell}.
$$
Hence, with probability at least $1-(c_3 {\cal L}\eps)^{c_4 \ell}$, there are at least $\ell/2$ indices $i \in I_{j(\ell)}$ such that
$$
|x_i a_i| \geq \eps \frac{a_1}{2^{j+1}},
$$
and in particular,
$$
\sum_{i \in I_{j(\ell)}} |x_i a_i|^p \geq c_5^p\eps^p \ell \frac{a_1^p}{2^{jp}} \geq c_6^p \eps^p \sum_{i \in I_{j(\ell)}} a_i^p.
$$
Set
$$
\phi^p(k)= \sum_{\{\ell \geq k, \Lambda_\ell \not = \emptyset\}} \sum_{i \in I_{j(\ell)}}  a_{i}^{p},
$$
note that
$$
(c_7\|a\|_p)^p \leq \phi^p(1) \leq \|a\|_p^p
$$
and that by the union bound, for every integer $k$, with probability at least $1-(c_8{\cal L}\eps)^{c_9 k}$,
$$
\left\|\sum_{i=1}^n a_i x_i e_i\right\|_p \geq \left(\sum_{\{\ell \geq k, \Lambda_\ell \not = \emptyset\}} \sum_{i \in I_{j(\ell)}} |a_i x_i|^p \right)^{1/p} \geq
c_6 \eps \phi(k).
$$
All that is left to show is that for a well chosen absolute constant $c$ and for $k = (c\|a\|_p/\|a\|_\infty)^{p}$ one has that $\phi(k) \gtrsim \|a\|_p^p$. To that end, and because $\phi^p(1) \geq c_7^p \|a\|_p^p$, the claim follows if
$$
\sum_{\ell \leq k} \sum_{i \in I_{j(\ell)}} a_i^p \leq \frac{c_7^p}{2} \|a\|_p^p.
$$
By the exponential decay of $\|P_{I_j} a\|_\infty$,
\begin{align*}
& \sum_{\{\ell < k, \Lambda_\ell \not = \emptyset\}} \sum_{i \in I_{j(\ell)}} a_i^p \leq \sum_{\{\ell < k, \Lambda_\ell \not = \emptyset\}} |I_{j(\ell)}| \|P_{I_{j(\ell)}}a\|_\infty^p \\
= & \sum_{\{\ell < k, \Lambda_\ell \not = \emptyset\}} |\ell| \|P_{I_{j(\ell)}}a\|_\infty^p \leq c_{10}^p k \|a\|_\infty^p
\end{align*}
from which the wanted estimate follows immediately for our choice of $k$.
\endproof

\section{Concluding Remarks}
Finally, let us describe how coordinate small-ball estimates \emph{should not be established}. Unfortunately, up to this point, the only known way of obtaining such estimates was this (suboptimal) way.

\subsection{The wrong way} \label{sec:wrong-way}
The standard way in which coordinate small-ball estimates have been established was based on the following simple observation. Consider a vector $x \in \R^n$ that satisfies $\|x\|_2 \geq \alpha \sqrt{n}$ for some $\alpha>0$ and for the sake of simplicity assume that $x_1 \geq x_2 \geq ... \geq x_n \geq 0$. Clearly, having any estimate on $\|x\|_2$ says nothing about the number of large coordinates that $x$ has; however, if the contribution to $\|x\|_2$ made by the $k=\beta n$ largest coordinates of $x$ is smaller than $\alpha \sqrt{n}/2$ then $x$ is ``well spread". Indeed, on the one hand
$$
\sum_{i=k+1}^n x_i^2 \geq \frac{3}{4}\alpha^2 n,
$$
and on the other
$$
x_k^2 \leq \frac{\alpha^2 n}{4k}.
$$
Therefore, by a Paley-Zygmund type argument, a proportional number ($\sim_{\alpha, \beta} n$) of the $x_i$'s are at least $c(\alpha,\beta)$.

 \vskip0.4cm

Obtaining a coordinate small-ball estimate in this way is particularly appealing in light of Theorem \ref{main-compl-sb}: because we know that $\|TX\|_2$ is likely to be large, it seems like half the job is already done. However, there are two crucial reasons why, despite the appeal, this is the wrong approach. Firstly, it gives no flexibility: one has no control on the proportion of nontrivial coordinates that the vector has, nor on the lower bound on the absolute values of these coordinates; in particular, there is no hope of proving Theorem \ref{th2-main-s2} using this type of argument. Secondly, while lower bounds are, in some sense, universal, upper bounds---which play an integral part in the argument and are based on tail estimates---are clearly not. In this case, given an orthonormal basis $(u_i)_{i=1}^m$ the necessary upper bound is on $\|(\inr{TX,u_i})_{i=1}^m\|_{[k]}$, where we set
$$
\|x\|_{[k]} = \max_{|I|=k} \left(\sum_{i \in I} x_i^2 \right)^{\frac{1}{2}}.
$$
Upper estimates of this kind hold with reasonable probability only in very special cases. In fact, even when $X$ is a Gaussian vector, the resulting probability estimate is weaker than, say, the one in Theorem \ref{th2-main-s2}; and for more heavy-tailed random vectors estimates on $\|(\inr{TX,u_i})\|_{[k]}$ are completely useless.

%
%

\vskip0.4cm
As a general principle,
\begin{framed}
It is wrong to try to establish coordinate small-ball estimates (which are lower bounds) using an argument that is based on ``large deviations".  Such a method may lead to nontrivial bounds only for very nice random vectors, and the bounds will be suboptimal even in those cases.
\end{framed}

\bibliographystyle{plain}
\bibliography{bib}

\begin{thebibliography}{10}

\bibitem{Ball1}
Keith Ball.
\newblock Logarithmically concave functions and sections of convex sets in
  {${\bf R}^n$}.
\newblock {\em Studia Math.}, 88(1):69--84, 1988.

\bibitem{Bou1}
J.~Bourgain.
\newblock On the distribution of polynomials on high-dimensional convex sets.
\newblock In {\em Geometric aspects of functional analysis (1989--90)}, volume
  1469 of {\em Lecture Notes in Math.}, pages 127--137. Springer, Berlin, 1991.

\bibitem{Bou2}
J.~Bourgain.
\newblock On the isotropy-constant problem for ``{PSI}-2''-bodies.
\newblock In {\em Geometric aspects of functional analysis}, volume 1807 of
  {\em Lecture Notes in Math.}, pages 114--121. Springer, Berlin, 2003.

\bibitem{BGVV}
Silouanos Brazitikos, Apostolos Giannopoulos, Petros Valettas, and
  Beatrice-Helen Vritsiou.
\newblock {\em Geometry of isotropic convex bodies}, volume 196 of {\em
  Mathematical Surveys and Monographs}.
\newblock American Mathematical Society, Providence, RI, 2014.

\bibitem{CFM}
D.~Cordero-Erausquin, M.~Fradelizi, and B.~Maurey.
\newblock The ({B}) conjecture for the {G}aussian measure of dilates of
  symmetric convex sets and related problems.
\newblock {\em J. Funct. Anal.}, 214(2):410--427, 2004.

\bibitem{FocRau}
Simon Foucart and Holger Rauhut.
\newblock {\em A mathematical introduction to compressive sensing}.
\newblock Applied and Numerical Harmonic Analysis. Birkh\"{a}user/Springer, New
  York, 2013.

\bibitem{GolLoa}
Gene~H. Golub and Charles~F. Van~Loan.
\newblock {\em Matrix computations}.
\newblock Johns Hopkins Studies in the Mathematical Sciences. Johns Hopkins
  University Press, Baltimore, MD, fourth edition, 2013.

\bibitem{Kl}
B.~Klartag.
\newblock On convex perturbations with a bounded isotropic constant.
\newblock {\em Geom. Funct. Anal.}, 16(6):1274--1290, 2006.

\bibitem{KV}
B.~Klartag and R.~Vershynin.
\newblock Small ball probability and {D}voretzky's theorem.
\newblock {\em Israel J. Math.}, 157:193--207, 2007.

\bibitem{LO}
Rafa\l Lata\l~a and Krzysztof Oleszkiewicz.
\newblock Small ball probability estimates in terms of widths.
\newblock {\em Studia Math.}, 169(3):305--314, 2005.

\bibitem{LMOT}
Rafal Latala, Piotr Mankiewicz, Krzysztof Oleszkiewicz, and Nicole
  Tomczak-Jaegermann.
\newblock Banach-{M}azur distances and projections on random subgaussian
  polytopes.
\newblock {\em Discrete Comput. Geom.}, 38(1):29--50, 2007.

\bibitem{LPP}
Galyna Livshyts, Grigoris Paouris, and Peter Pivovarov.
\newblock On sharp bounds for marginal densities of product measures.
\newblock {\em Israel J. Math.}, 216(2):877--889, 2016.

\bibitem{LM3}
G.~Lugosi and S.~Mendelson.
\newblock Near-optimal mean estimators with respect to general norms.
\newblock {\em Probab. Theory Relat. Fields}, to appear.

\bibitem{LM1}
G.~Lugosi and S.~Mendelson.
\newblock Risk minimization by median-of-means tournaments.
\newblock {\em J. Eur. Math. Soc.}, to appear.

\bibitem{LM2}
G\'{a}bor Lugosi and Shahar Mendelson.
\newblock Sub-{G}aussian estimators of the mean of a random vector.
\newblock {\em Ann. Statist.}, 47(2):783--794, 2019.

\bibitem{M2}
S.~Mendelson.
\newblock Approximating the covariance ellipsoid.
\newblock {\em Preprint}.

\bibitem{M1}
S.~Mendelson.
\newblock An optimal unrestricted learning procedure.
\newblock {\em Preprint}.

\bibitem{MZ}
S.~Mendelson and N.~Zhivotovskiy.
\newblock Robust covariance estimation under \mbox{$L_4-L_2$} norm equivalence.
\newblock {\em Preprint}.

\bibitem{Men}
Shahar Mendelson.
\newblock Learning without concentration.
\newblock {\em J. ACM}, 62(3):Art. 21, 25, 2015.

\bibitem{MMP}
Shahar Mendelson, Emanuel Milman, and Grigoris Paouris.
\newblock Generalized dual {S}udakov minoration via dimension-reduction---a
  program.
\newblock {\em Studia Math.}, 244(2):159--202, 2019.

\bibitem{MRW}
Shahar Mendelson, Holger Rauhut, and Rachel Ward.
\newblock Improved bounds for sparse recovery from subsampled random
  convolutions.
\newblock {\em Ann. Appl. Probab.}, 28(6):3491--3527, 2018.

\bibitem{NY}
Assaf Naor and Pierre Youssef.
\newblock Restricted invertibility revisited.
\newblock In {\em A journey through discrete mathematics}, pages 657--691.
  Springer, Cham, 2017.

\bibitem{Paouris12}
Grigoris Paouris.
\newblock On the isotropic constant of marginals.
\newblock {\em Studia Math.}, 212(3):219--236, 2012.

\bibitem{P-sb}
Grigoris Paouris.
\newblock Small ball probability estimates for log-concave measures.
\newblock {\em Trans. Amer. Math. Soc.}, 364(1):287--308, 2012.

\bibitem{PP-IMA}
Grigoris Paouris and Peter Pivovarov.
\newblock Randomized isoperimetric inequalities.
\newblock In {\em Convexity and concentration}, volume 161 of {\em IMA Vol.
  Math. Appl.}, pages 391--425. Springer, New York, 2017.

\bibitem{Rom}
Justin Romberg.
\newblock Compressive sensing by random convolution.
\newblock {\em SIAM J. Imaging Sci.}, 2(4):1098--1128, 2009.

\bibitem{RV}
Mark Rudelson and Roman Vershynin.
\newblock Small ball probabilities for linear images of high-dimensional
  distributions.
\newblock {\em Int. Math. Res. Not. IMRN}, (19):9594--9617, 2015.

\end{thebibliography}

\newpage

\appendix

\section{Examples of vectors that satisfy the \it{SBA}} \label{app:log-concave0}
Here we give examples of several generic random vectors that satisfy the {\it SBA}. This is far from being an exhaustive list and should be viewed only as an indication to the fact that the {\it SBA} is a property shared by many natural random vectors.

\begin{description}
\item{$(1)$} Let $X= (\xi_{1}, \cdots, \xi_{n})$ where the $ \xi_{i}$'s are independent random variables with densities bounded by $ {\cal{L}}$. It was shown in \cite{RV} that $X$ satisfies the {\it SBA} with constant $ c{\cal{L}}$, where $ c>0$ is an absolute constant. 
    
    This fact was further extended in \cite{LPP,PP-IMA}; most notably, it was shown in \cite{PP-IMA} that if the coordinates of $X=(\xi_i)_{i=1}^n$ are independent random variables with densities bounded by $1$ and the coordinates of $Y=(\eta_i)_{i=1}^n$ are uniformly distributed in $[-\frac{1}{2} , \frac{1}{2}]$, then for every semi-norm $ \| \cdot \|$ and $ t>0$,
\begin{equation} \label{SB-1-1}
\mathbb P ( \|X\|\leq t ) \leq \mathbb P ( \| Y \|\leq t ).
\end{equation}
In particular, among all such vectors the `worse' small-ball behaviour---with respect to \emph{any} semi-norm---is exhibited by the uniform measure on the cube $[-\frac{1}{2} , \frac{1}{2}]^n$.

Observe that for the Euclidean norm, the small-ball behaviour of $Y$ and of the standard Gaussian vector $G$ is the same up to absolute constants.

\item{$(2)$} \emph{Perturbations:} It is standard to verify that if $X$ satisfies the {\it SBA} with a constant ${\cal L}$ and $W$ is an arbitrary random vector that is independent of $X$, then $ W+\delta X$ satisfies {\it SBA} with a constant depending on $ \delta$ and ${\cal L}$.

\item{$(3)$} The question of whether there is a constant ${\cal L}$ such that \emph{any} isotropic log-concave random vector satisfies the {\it SBA} with constant ${\cal L}$ is equivalent to Bourgain's celebrated Hyperplane Conjecture (see \cite{Bou1} and the discussion in \cite{Kl} and \cite{BGVV})).

Thanks to the extensive study of log-concave measures and the connection the SBA has with the Hyperplane conjecture for such measures, there are some important examples of isotropic, log-concave random vectors that are known to satisfy the {\it SBA} with an absolute constant:

\noindent $\bullet$ If $X$ is also $1$-unconditional  (see \cite{MMP}, section 8.2);

\noindent $\bullet$ If $X$ is also subgaussian (\cite{Bou1,Bou2});

\noindent $\bullet$ If $X$ is also supergaussian (this follows from results of \cite{Paouris12}).
\end{description}

\section{Proof of Remark \ref{rem:sharp}}  \label{app:log-concave}

The proof requires some additional notation. Let $X$ be a random vector in $ \R^{n}$ and let $p\geq 1$. The $Z_{p}$ body of $X$ is defined as the (centrally-symmetric) convex body whose support function is
\begin{equation} \label{Zp-1}
h_{Z_{p} (X) } ( \theta) = \left( \E|\inr{X, \theta}|^{p} \right)^{\frac{1}{p}} , \ \ \ \theta \in S^{n-1} .
\end{equation}
It is straightforward to verify that if $ T :\R^{n} \to \R^{m}$ is a linear operator then
\begin{equation} \label{Zp-2}
Z_{p} (T X) = T Z_{p} ( X).
\end{equation}

\begin{Lemma} \label{lemma:A1}
There are absolute constants $c_1$ and $c_2$ for which the following holds.
Let $ X$ be a centred log-concave random vector in $\R^{n}$ that satisfies the {\it SBA} with constant $ {\cal{L}}$. For any $ T\in GL_{n}$ and $ F\in {\cal G}_{n,k}$ one has
\begin{equation} \label{log-conc-2}
\frac{c_{1}}{ | {\rm det}[(P_{F} T) ( P_{F} T)^{\ast}] |^{\frac{1}{2k}} }   \leq   f_{P_{F}TX}^{\frac{1}{k}}(0)  \leq \frac{c_{2}{\cal{L}} }{ (  {\rm det}[(P_{F} T) ( P_{F} T)^{\ast}] )^{\frac{1}{2k}} },
\end{equation}
where the left-hand side holds true under the additional assumption that $X$ is isotropic.
\end{Lemma}

The proof of Lemma \ref{lemma:A1} is based on two facts. The first is a standard observation from linear algebra:
let $T: \R^{n} \to \R^{k}$, set $E= {\rm ker}(T)^{\perp}= {\rm im} (T^{\ast} )  $ and denote by $ T|_{E}$ the restriction of $T$ to $E$. Then for any compact set $K \subset \R^n$,
\begin{equation} \label{Lin00}
{\rm vol}(T K) = {\rm det} ( T T^{\ast} )\cdot {\rm vol}(P_{E} K) .
\end{equation}
The second observation is Proposition~3.7 from \cite{P-sb}: If $X$ is a centred, log-concave random vector then
\begin{equation} \label{Zp-3}
f_{X}^{\frac{1}{n}} (0)\sim {\rm vol}^{-\frac{1}{n}}(Z_{n} ( X) ).
\end{equation}

\noindent {\bf Proof of Lemma \ref{lemma:A1}.}
By the Prekop\'a-Leindler inequality, for every linear operator $S$, the random vector $SX$ is also log-concave and centred. Hence,
using \eqref{Zp-2}, \eqref{Zp-3} and \eqref{Lin00},
 it is evident that
\begin{align*}
& f_{P_{F}TX}^{\frac{1}{k}}(0) \sim  ({\rm vol} (Z_{k} ( P_{F} T X) ))^{-\frac{1}{k}} \sim  ({\rm vol} (P_{F} T Z_{k} (  X) ))^{-\frac{1}{k}}
\\
\sim & \frac{1}{ (   {\rm det}[(P_{F} T) ( P_{F} T)^{\ast}] )^{\frac{1}{2k}} }   ({\rm vol} (  P_{E}  Z_{k} (  X) ))^{-\frac{1}{k}}  \sim
   \frac{1}{ (  {\rm det}[(P_{F} T) ( P_{F} T)^{\ast}])^{\frac{1}{2k}} }   ({\rm vol} (Z_{k} ( P_{E}  X) ))^{-\frac{1}{k}}
\\
\sim &  \frac{ f_{P_{E}X}^{\frac{1}{k}}(0)}{ (   {\rm det}[(P_{F} T) ( P_{F} T)^{\ast}] )^{\frac{1}{2k}} }.
\end{align*}
Clearly, $ f_{P_{E}X}^{\frac{1}{k}}(0) \leq {\cal{L}}$, which proved the right-hand side inequality in \eqref{log-conc-2}. Moreover if $X$ is an isotropic log-concave random vector in $\R^{n}$ then $f_{X}^{\frac{1}{n}} (0 ) \geq c $, where $c$ is an absolute constant (see, e.g. \cite{Ball1}). And since $ P_{F} X$ is also isotropic when $X$ is, the left-hand side inequality in \eqref{log-conc-2} follows.
\endproof

Combining \eqref{log-conc-2} and \eqref{P-sb-1} it is evident that:
\begin{Proposition}
There are absolute constants $c_1$ and $c_2$ for which the following holds.
Let $X$ be an isotropic log-concave random vector in $\R^{n}$ that satisfies the {\it SBA} with constant $ {\cal{L}}$ and let $ T: \R^{n} \to \R^{m}$ be a linear operator. Then
\begin{equation} \label{log-conc}
\left( \mathbb E \| T G/ (c_{2} {\cal{L}}) \|_{2}^{-k} \right)^{-\frac{1}{k}} \leq \left( \mathbb E \| T X \|_{2}^{-k} \right)^{-\frac{1}{k}} \leq \left( \mathbb E \| T G/c_{1} \|_{2}^{-k} \right)^{-\frac{1}{k}}.
\end{equation}
\end{Proposition}


\end{document}